%% file: fk.tex
\documentclass[11pt,reqno,oneside]{amsart}

\usepackage{graphicx}    % standard LaTeX graphics tool when including figure files
\usepackage{amsfonts}
\usepackage{verbatim}
\usepackage{amssymb}
\usepackage{amsmath}
\usepackage{hyperref}
\usepackage[numbers,sort&compress]{natbib}

\usepackage[it,small]{caption}  % fonts for captions in figures and tables

\newcommand{\eps}{\varepsilon}
\newcommand{\D}{\:\mathrm{d}}

\newcommand{\Z}{{\mathbb Z}}
\newcommand{\R}{{\mathbb R}}

\newcommand{\Oh}{{\mathcal O}}

\newcommand{\E}{\mathcal{E}}

\DeclareMathOperator*{\argmin}{arg\,min}
\DeclareMathOperator*{\sumprime}{{\sum}^\prime}

\newcommand{\half}{{\textstyle \frac{1}{2}}}
\newcommand{\fourth}{{\textstyle \frac{1}{4}}}

\renewcommand{\vec}[1]{\mathbf{#1}}

\newtheorem{lemma}{Lemma}[section]
\newtheorem{theorem}{Theorem}[section]

\numberwithin{equation}{section}

\allowdisplaybreaks[4]

\begin{document}

\title[Adaptivity for the Quasicontinuum Approximation of a Frenkel-Kontorova Model]
{Error Estimation and Atomistic-Continuum Adaptivity\\
for the Quasicontinuum Approximation\\
of a Frenkel-Kontorova Model}
\author{Marcel Arndt}
\author{Mitchell Luskin}
\begin{abstract}
  We propose and analyze a goal-oriented {\em a posteriori} error estimator for the
  atomistic-continuum modeling error
  in the quasicontinuum method.  Based on this error estimator, we develop an
  algorithm which adaptively determines the atomistic and continuum regions
  to compute a quantity of interest to within a given tolerance.
  We apply the algorithm to the computation of the structure of a
  crystallographic defect described by a Frenkel-Kontorova model and present the
  results of numerical experiments.  The numerical results show that our method
  gives an efficient estimate of the error and a nearly optimal
  atomistic-continuum modeling strategy.
\end{abstract}
\maketitle

\input{intro}        % Introduction
\input{qcapprox}     % Quasicontinuum Approximation
\input{fkmodel}      % Frenkel-Kontorova: Atomistic and Quasicontinuum Model
\input{problem}      % Primal and Dual Problems
\input{errorcont}    % Error due to Continuum Formulation
\input{numericscont} % Numerical Results for Atomistic/Continuum Modelling
\appendix
%\input{linear}      % Linearization
\input{matrix}      % Matrix Definitions

\bibliographystyle{hsiam}
\bibliography{../Literature/marcel}

\end{document}

%% file: intro.tex
% $Id: intro.tex,v 1.16 2007/04/11 19:50:06 arndt Exp $

%%%%%%%%%%%%%%%%%%%%%%%%%%%%%%%%%%%%%%%%%%%%%%%%%%%%%%%%%%%%%%%%%%%%%%%%%%%%%%%%
%%%%%%%%%%%%%%%%%%%%%%%%%%%%%%%%%%%%%%%%%%%%%%%%%%%%%%%%%%%%%%%%%%%%%%%%%%%%%%%%
%%%%%%%%%%%%%%%%%%%%%%%%%%%%%%%%%%%%%%%%%%%%%%%%%%%%%%%%%%%%%%%%%%%%%%%%%%%%%%%%
\section{Introduction}

The quasicontinuum (QC) method \cite{TadmorMillerPhillipsOrtiz:1999,
  TadmorOrtizPhillips:1996, TadmorPhilipsOrtiz:1996} has been successfully used
to efficiently couple atomistic and continuum models for crystalline solids and
offers the possibility of computing mesoscale or macroscale properties by a
nearly minimal number of degrees of freedom.  Accurate modeling requires that an
atomistic model be used in regions with highly non-uniform deformations such as
around dislocations, whereas a continuum model can be used in regions with
nearly uniform deformations to reduce the number of degrees of freedom.

It is usually not known {\em a priori} which regions of some specimen undergo
uniform deformations and which do not, so {\em a posteriori} error estimation
is important for the design of efficient numerical approximations
by the quasicontinuum method.
Since the purpose of a computation is often to obtain the
value of a (usually local) quantity of interest to a desired
error tolerance rather than to obtain a solution to a desired
error tolerance for a global norm, there has been great interest in
the development of goal-oriented error estimators for many problems.
They are based on duality techniques and have been developed and used to
adaptively refine finite element approximations of continuum
problems~\cite{AinsworthOden:2000, BangerthRannacher:2003} and to study and
control modeling error~\cite{OdenVemaganti:2000}.

In this paper, we extend this approach to develop an {\em a posteriori} error
estimator for the quasicontinuum method which quantifies the atomistic-continuum modeling error for
a goal function and allows for an adaptive decision about which regions can
be accurately modeled as a continuum and which regions need to be modeled atomistically.
Methods to determine the optimal mesh size within the continuum region will be
studied in a forthcoming paper.

Crystallographic defects~\cite{ChaikinLubensky:2000} provide a challenge
to validate atomistic-continuum error estimators and adaptivity.  No
such error estimators and adaptive methods currently
exist for fully three-dimensional crystals.  As a step in
this direction, we develop a rigorous theory for a simple one-dimensional
atomistic model for a defect that is a modification of the Frenkel-Kontorova
model~\cite{Marder:2000}.  We add next-nearest-neighbor harmonic
interactions between the atoms to the nearest-neighbor harmonic interactions
between the atoms in the classical Frenkel-Kontorova model.

{\em A priori} analyses for various quasicontinuum approximations have
been given in \cite{Lin:2003, Lin:2005, DobsonLuskin:2006, BlancLeBrisLegoll:2005,
EMing:2004, ELuYang:2006, OrtnerSueli:2006a, KnapOrtiz:2001}.
An {\em a posteriori} analysis for a slightly different one-dimensional quasicontinuum
approximation is given in \cite{OrtnerSueli:2006b}.
The development and application of a goal-oriented error estimator for mesh coarsening in
a two-dimensional quasicontinuum method
is reported in \cite{OdenPrudhommeRomkesBauman:2006, OdenPrudhommeBauman:2006}.

Let us mention that the continuum model used in the QC method, which coincides
with the model obtained by the classical thermodynamic limit, is by far not
the only reasonable continuum model to use.  A method to derive continuum models
which approximate atomistic models up to an arbitrarily high order has been
proposed in \cite{ArndtGriebel:2005}.

The paper is organized as follows.  In Section~\ref{SecQCApprox}, we give a
general formulation of the one-dimensional quasicontinuum
approximation~\cite{TadmorOrtizPhillips:1996} that includes not only two-body
and three-body
potentials, but also many body potentials such as the embedded atom
potential~\cite{DawBaskes:1983, DawBaskes:1984}.  In Section~\ref{SecFKModel},
we describe our extension of the Frenkel-Kontorova model and its quasicontinuum
approximation.  In Section~\ref{SecProblem}, we introduce the primal and dual
problems for our model and formulate our approach to goal-oriented error
estimation.

Next, in Section~\ref{SecLocNonloc} we extend the approach in~\cite{OdenPrudhomme:2002} to
develop an error estimator for atomistic-continuum modeling.
This first error estimator does not allow a decomposition
among the atoms that can be used for atomistic-continuum adaptivity,
so we propose and analyze a less accurate second error estimator that
does allow such a decomposition.

Finally, in Section~\ref{SecNumericalResults} we propose an adaptive atomistic-continuum
modeling algorithm and show that it gives an efficient estimate of
the modeling error and a nearly optimal atomistic-continuum
modeling strategy for the computation of defect structure.

%% file: qcapprox.tex
% $Id: qcapprox.tex,v 1.22 2007/04/11 19:50:06 arndt Exp $

%%%%%%%%%%%%%%%%%%%%%%%%%%%%%%%%%%%%%%%%%%%%%%%%%%%%%%%%%%%%%%%%%%%%%%%%%%%%%%%%
%%%%%%%%%%%%%%%%%%%%%%%%%%%%%%%%%%%%%%%%%%%%%%%%%%%%%%%%%%%%%%%%%%%%%%%%%%%%%%%%
%%%%%%%%%%%%%%%%%%%%%%%%%%%%%%%%%%%%%%%%%%%%%%%%%%%%%%%%%%%%%%%%%%%%%%%%%%%%%%%%
\section{Quasicontinuum Approximation}  \label{SecQCApprox}

The departure point for the QC approximation is the potential energy of the
atomistic system.  The potential energy that is utilized
fully models the properties of the
system. The local minima of the potential energy model the metastable states of
the system, and the potential energy can be used in Newton's equations of motion
to model the dynamical behavior.

The QC method approximates the potential energy of the atomistic system in two
steps.  First, we develop a continuum potential energy that will be used in the
adaptively determined continuum region,
and we then show how to
reduce the degrees of freedom in the continuum region.

%%%%%%%%%%%%%%%%%%%%%%%%%%%%%%%%%%%%%%%%%%%%%%%%%%%%%%%%%%%%%%%%%%%%%%%%%%%%%%%%
\subsection{The Atomistic System}

We assume that the atomistic system has $2M$
atoms with deformation given by $\vec
y^a= (y^a_{-M+1},\dots y^a_M)\in\R^{2M}$.
Without loss of generality, we assume that the atoms are
ordered so that their positions satisfy $y^a_i < y^a_{i+1}$.
Furthermore, we assume that the atomistic total
potential energy, $\E^a(\vec y^a),$ can be written as a sum over
potential energies associated with each atom, $\E^a_i(\vec y^a),$
so that
\begin{align}\label{decomp}
  \E^a(\vec y^a) = \sum_{i=-M+1}^M \E^a_i(\vec y^a).
\end{align}
This decomposition can be found for most empirical potentials,
including embedded atom potential energies~\cite{DawBaskes:1983,
DawBaskes:1984}. For example, if the
atomistic total
potential energy $\E^a(\vec y^a)$ is given by
\begin{equation}\label{class}
\E^a(\vec y^a)=\sum_{i<j}\psi(y^a_j-y^a_i),
\end{equation}
where $\psi(r)$ is an empirical two-body potential
energy, then we can obtain the decomposition
\eqref{decomp} by taking
\begin{equation}\label{at}
\E^a_i(\vec y^a)=\frac 12 \sum_{j\ne i}
\psi(y^a_j-y^a_i).
\end{equation}
We note that $\E^a_i(\vec y^a)$ can also contain contributions from
external forces, such as for the Frenkel-Kontorova model described in
Section~\ref{SecFKModel}, and can thus depend on $i.$

%%%%%%%%%%%%%%%%%%%%%%%%%%%%%%%%%%%%%%%%%%%%%%%%%%%%%%%%%%%%%%%%%%%%%%%%%%%%%%%%
\subsection{The Atomistic-Continuum Energy}
For any deformation $\vec
y^a\in\R^{2M},$ we let $L^{i,i+1}\vec y^a \in \R^\Z$ denote the
linear extrapolation of the atomistic positions $y^a_i$ and
$y^a_{i+1}$ given by
\begin{equation} \label{EqDefL}
  (L^{i,i+1} \vec y^a)_k = (k-i) y^a_{i+1} + (i+1-k) y^a_i\qquad\text{for }k=-\infty,\dots,\infty.
\end{equation}
The continuum potential energy $\E^c_i(\vec y^a)$ of atom $i$ is obtained from the
average of the atomistic potential energy $\E^a_i$ evaluated at the
extrapolations $L^{i-1,i} \vec y^a$ and $L^{i,i+1} \vec y^a$ by
\begin{align} \label{EqDefEaci}
  \E^c_i(\vec y^a) :=   \half \E^a_i(L^{i-1,i} \vec y^a)
                      + \half \E^a_i(L^{i,i+1} \vec y^a),
\end{align}
where we note that the domain of $\E^a_i$ has been expanded to
the infinite periodic atomistic systems in the range of $L^{i-1,i}$
and $L^{i,i+1}.$  We assume that $\E^a_i$ is finite for
infinite periodic atomistic systems, which is true for
\eqref{at} when the two-body potential $\psi(r)$ decays fast
enough so that $\sum_{k=1}^\infty \psi(kr)$ is finite for $r\ne 0.$
At the endpoints of the chain, the extrapolation can be done only to one side,
so we neglect the undefined part and define
\begin{align}
  \E^c_{-M+1}(\vec y^a) := \half \E^a_{-M+1}(L^{-M+1,-M+2} \vec y^a)
  \qquad \text{and} \qquad
  \E^c_M(\vec y^a) := \half \E^a_M(L^{M-1,M} \vec y^a).
\end{align}

We then decide for each atom $i$ whether to model its energy
atomistically by $\E^a_i(\vec y^a)$ or as a continuum by
$\E^c_i(\vec y^a)$.
We thus obtain for the whole chain the atomistic-continuum energy
\begin{equation}\label{EqEcDef}
\begin{split}
  \E^{ac}(\vec y^a)
  := & \sum_{i=-M+1}^M \delta^a_i \E^a_i(\vec y^a)
     + \sum_{i=-M+1}^M \delta^c_i \E^c_i(\vec y^a) \\
  =  & \sum_{i=-M+1}^M \delta^a_i \E^a_i(\vec y^a)
     + \half \sum_{i=-M+2}^M     \delta^c_i \E^a_i(L^{i-1,i} \vec y^a)
     + \half \sum_{i=-M+1}^{M-1} \delta^c_i \E^a_i(L^{i,i+1} \vec y^a),
     \end{split}
\end{equation}
where
\begin{align}
\delta_i^a = \begin{cases}
  1 & \text{if atom $i$ is modeled atomistically,} \\
  0 & \text{if atom $i$ is modeled as continuum,}
\end{cases}
\qquad \text{and} \qquad \delta_i^c = 1-\delta_i^a.
\end{align}
This approximation allows for a slightly faster evaluation of the
energy and its derivatives, especially if $\E^a_i$ is long-ranged.
However, it reveals its full strength only after the
quasicontinuum coarsening to be described next.
We note that sometimes atomistic degrees of freedom and energies
are referred to as nonlocal and continuum degrees of freedom and energies
are referred to as local \cite{TadmorOrtizPhillips:1996}.

\subsection{Repatoms: Reduction of Degrees of Freedom}
The quasicontinuum method allows a reduction of the number of degrees
of freedom in the continuum region. To this end,
we choose so-called {\em representative atoms}, or more
briefly called {\em repatoms}. The repatoms are a subset of
the original atoms. The quasicontinuum approximation of the
energy is defined completely in terms of the repatoms.

We choose the repatoms by defining indices $\ell_j$ for $j=-N+1,
\ldots, N$ where
\[ -M+1=\ell_{-N+1} < \cdots < \ell_j <
\ell_{j+1} < \cdots < \ell_N=M.
\]
  The atoms at $y^a_i$ for $i=\ell_{-N+1},
\ell_{-N+2}, \ldots, \ell_N$ are repatoms, and all of the remaining atoms are
non-repatoms.   We have that
\begin{align}
  \nu_j = \ell_{j+1}-\ell_j
\end{align}
gives the number of atomistic intervals between the repatoms $\ell_j$ and
$\ell_{j+1}.$ We require that the
chain not be coarsened in the atomistic regions, which precisely means that
$\delta^c_{\ell_j} = \delta^c_{\ell_j+1} = \ldots = \delta^c_{\ell_{j+1}}
= 1$ whenever $\nu_j>1$.

Finally, the interactions of the atomistic energy only partially reach into the continuum
part if the atomistic potential has a finite cutoff radius. To allow for
an exact calculation of this energy without atomistic
interpolation, we require that these regions are not coarsened as well.  As we
will see in the next subsection, the atomistic next-nearest-neighbor interactions from the
Frenkel-Kontorova model studied in this paper reach two atoms into the
continuum part.  Hence, we require that $\nu_{j-2} = \nu_{j-1} = \nu_j = \nu_{j+1} =
1$ whenever $\delta^a_{\ell_j}=1$. Other potential energies in general require
similar conditions that depend on their cut-off radius.

We denote the position of the $j$-th repatom by
$y^{qc}_j=y^a_{\ell_j}$ and the vector of all repatoms by $\vec
y^{qc}\in\R^{2N}$.

%%%%%%%%%%%%%%%%%%%%%%%%%%%%%%%%%%%%%%%%%%%%%%%%%%%%%%%%%%%%%%%%%%%%%%%%%%%%%%%%
\subsection{The Quasicontinuum Energy}

Now we define the quasicontinuum energy. To this end, the missing non-repatoms
are implicitly reconstructed. We will see later that this helps to set up the QC
model, but needs not be done for the actual computation.

The reconstruction is done by a linear interpolation
between the nearest repatom to the right and to the left. That is,
the vector of all atomistic positions is computed from the vector
$\vec y^{qc}$ of repatom positions by
\begin{align}
  I: \R^{2N} \to \R^{2M}, \qquad
  (I \vec y^{qc})_{\ell_j+m} := \frac{\nu_j-m}{\nu_j} y^{qc}_j
                            + \frac{m      }{\nu_j} y^{qc}_{j+1},
  \qquad m=0,\ldots,\nu_j.
\end{align}
We note that
\begin{align}
  y^{qc}_j = (I \vec y^{qc})_{\ell_j}.
\end{align}

The underlying idea is that in regions where the lattice spacing
of the atoms is nearly constant, this interpolation is very close
to the actual atomistic positions and therefore leads to a good
approximation of the total energy. Only a few repatoms are needed
in these regions. This exactly corresponds to mesh coarsening in
classical finite element approximations of continuum
models.  On the other hand, in regions where the lattice spacing
is non-uniform, such as around a dislocation, all atoms must be chosen
to be repatoms to obtain sufficient accuracy.
This guarantees that the full resolution of the
atomistic model in the critical regions is retained and
corresponds to a high refinement in classical finite element
continuum models.

% This is commonly referred to as the Cauchy-Born hypothesis.

%Altogether, the restriction to repatoms
%reduces the number of degrees of freedom to speed up the solution process while
%retaining the high resolution in those regions which are necessary to obtain an
%accurate solution.

We define the QC approximation of the total energy to be
\begin{align} \label{EqQC1}
  \E^{qc}(\vec y^{qc}) := \E^{ac}(I \vec y^{qc}).
\end{align}
Now \eqref{EqQC1} has to be reformulated such that it can be
computed efficiently, without the overhead of the interpolation.
Most atomistic potentials are invariant to translations, a
property that allows us to simplify \eqref{EqQC1} considerably.
For any translationally invariant energy $\E^a_i$, we have that
$\E^a_i(L^{i,i+1} \vec y^a) = \phi_i(y^a_{i+1}-y^a_{i})$ and
$\E^a_i(L^{i-1,i} \vec y^a) = \phi_i(y^a_{i}-y^a_{i-1})$ for some
function $\phi_i$. If these functions
$\phi_i$ coincide, that is, $\phi_i=\phi_j$ for all
$i$ and $j,$ 
% Translational invariance in our context means that $\E^a_i(\vec
% y^a) = \E^a_i(\vec y^a+c)$ for any constant $c\in\R$. Thus,
% $\E^a_i(L^{i-1,i} \vec y^a)$ and $\E^a_i(L^{i,i+1} \vec y^a)$
% depend only on the quantities $y^a_i-y^a_{i-1}$ and
% $y^a_{i+1}-y^a_i$, respectively. Moreover, we assume that $\E^a_i
% = \E^a_j$ for all $i,j$. THE PREVIOUS SENTENCE IS CONFUSING. Thus,
we can write
\begin{align} \label{EqEnergyDensity}
  \E^a_i(L^{i-1,i} \vec y^a) = \phi \left( y^a_{i}-y^a_{i-1} \right)
  \qquad \text{and} \qquad
  \E^a_i(L^{i,i+1} \vec y^a) = \phi \left( y^a_{i+1}-y^a_i \right)
\end{align}
for some function $\phi:\R\to\R$. Here $\phi$ plays the role of a
continuum energy density and is given for the two-body potential
\eqref{class} by
\[
\phi(r)=\sum_{k=1}^\infty \psi(kr).
\]

Equations \eqref{EqEcDef}, \eqref{EqQC1}, and \eqref{EqEnergyDensity} lead to
\begin{align} \label{EqQC2}
  \E^{qc}(\vec y^{qc})
  & = \sum_{i=-M+1}^M \delta^a_i \E^a_i(I \vec y^{qc})
     + \half \sum_{i=-M+2}^M     \delta^c_i \E^a_i(L^{i-1,i} I \vec y^{qc})
     + \half \sum_{i=-M+1}^{M-1} \delta^c_i \E^a_i(L^{i,i+1} I \vec y^{qc}) \nonumber\\
  & = \sum_{i=-M+1}^M \delta^a_i \E^a_i(I \vec y^{qc})
     + \half \sum_{i=-M+2}^M \delta^c_i
       \phi( (I \vec y^{qc})_i    - (I \vec y^{qc})_{i-1} )\\
   &\qquad\qquad  + \half \sum_{i=-M+1}^{M-1} \delta^c_i
       \phi( (I \vec y^{qc})_{i+1} - (I \vec y^{qc})_i ).\nonumber
\end{align}
Because $I \vec y^{qc}$ is the linear interpolation between two repatoms $y^{qc}_j$ and
$y^{qc}_{j+1},$ we have
\begin{align}
  (I \vec y^{qc})_{i+1} - (I \vec y^{qc})_i = \frac{y^{qc}_{j+1}-y^{qc}_j}{\nu_j},
  \qquad i=\ell_j, \ldots, \ell_{j+1}-1.
\end{align}
Hence,
\begin{align} \label{EqQC3}
  \E^{qc}(\vec y^{qc}) =
  \sum_{i=-M+1}^M \delta^a_i \E^a_i(I \vec y^{qc}) +
  \sum_{j=-N+1}^{N-1} \omega_j \phi \left( \frac{y^{qc}_{j+1}-y^{qc}_j}{\nu_j} \right)
\end{align}
with weight factors
\begin{align}
  \omega_j = \half \nu_j \big( \delta^c_{\ell_j} + \delta^c_{\ell_{j+1}} \big)
  = \begin{cases}
    0     & \text{if both } y^{qc}_j \text{ and } y^{qc}_{j+1} \text{ are atomistic}, \\
    \half & \text{if exactly one of } y^{qc}_j \text{ and } y^{qc}_{j+1} \text{ is continuum}, \\
    \nu_j & \text{if both } y^{qc}_j \text{ and } y^{qc}_{j+1} \text{ are continuum}. \\
  \end{cases}
\end{align}
The first sum corresponds to the atomistic region which will
be a small region and is thus computationally inexpensive.
The second sum only involves at most $2N$ terms which is a considerable
reduction when $N\ll M.$

Note that the second term in formula \eqref{EqQC3} coincides with an integral
over the energy density $\phi$ as it occurs in finite element discretizations of
classical continuum mechanical models.  Hence the apparently unmotivated
definitions \eqref{EqEcDef} and \eqref{EqDefEaci} of the continuum energy here
result in what is commonly understood as a continuum energy.  The linear
interpolation operator $I$ resembles the Cauchy-Born hypothesis.

%% file: fkmodel.tex
% $Id: fkmodel.tex,v 1.23 2007/04/11 19:50:07 arndt Exp $

%%%%%%%%%%%%%%%%%%%%%%%%%%%%%%%%%%%%%%%%%%%%%%%%%%%%%%%%%%%%%%%%%%%%%%%%%%%%%%%%
%%%%%%%%%%%%%%%%%%%%%%%%%%%%%%%%%%%%%%%%%%%%%%%%%%%%%%%%%%%%%%%%%%%%%%%%%%%%%%%%
%%%%%%%%%%%%%%%%%%%%%%%%%%%%%%%%%%%%%%%%%%%%%%%%%%%%%%%%%%%%%%%%%%%%%%%%%%%%%%%%
\section{Frenkel-Kontorova Model} \label{SecFKModel}

\begin{figure}
\includegraphics[width=0.6\textwidth]{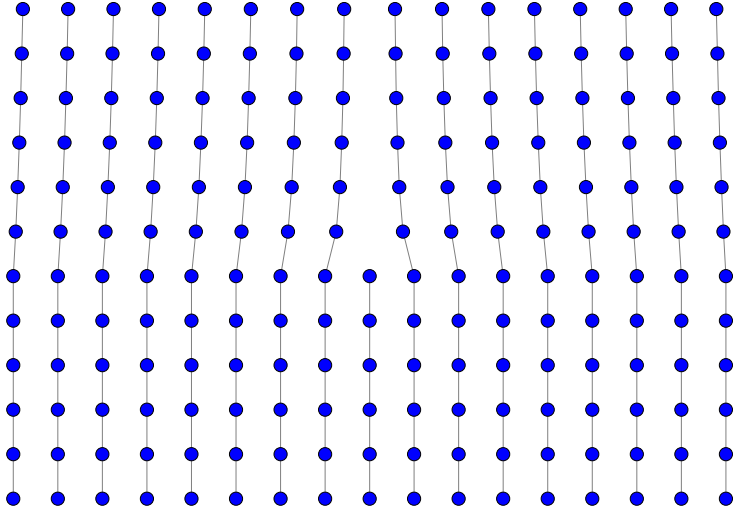}
\caption{Cross-section through a dislocation in a three-dimensional cubic
  lattice. The displayed plane repeats periodically in the three-dimensional
  crystal. Vertical bonds are shown by lines to emphasize the topological
  defect.}
\label{FigDisloc2D}
\end{figure}

Dislocations are lines in crystals which represent a defect
in the lattice structure~\cite{Marder:2000}, see Figure~\ref{FigDisloc2D}.
Typically, there is a core of small radius surrounding
the dislocation line where the lattice structure is highly deformed,
but the lattice structure is nearly uniform outside the core.
A simple one-dimensional model for a defect such as
a dislocation is given by
the Frenkel-Kontorova model~\cite{Marder:2000}.
Here, the elastic energy is modeled by harmonic interactions
between the atoms in the one-dimensional chain and the misfit
energy of the slip plane is modeled by a periodic potential.
A more accurate model of the same form is given by the
Peierls-Nabarro model~\cite{Kaxiras:2003}.

%%%%%%%%%%%%%%%%%%%%%%%%%%%%%%%%%%%%%%%%%%%%%%%%%%%%%%%%%%%%%%%%%%%%%%%%%%%%%%%%
\subsection{Atomistic Frenkel-Kontorova Model} \label{SecFKAtomistic}

\begin{comment}
We first model the core structure of a dislocation, which is a
one-dimensional singularity (with core) in three-dimensional physical space (see
Figure~\ref{FigDisloc2D}),
by a point singularity (with core) in a one-dimensional chain using the Frenkel-Kontorova
approximation~\cite{Marder:2000}. The chain models one layer of
atoms in a crystal that relaxes its energy above a rigid layer
of the crystal below, see Figure~\ref{FigFKModel}.
\end{comment}

\begin{figure}
\includegraphics[width=0.9\textwidth]{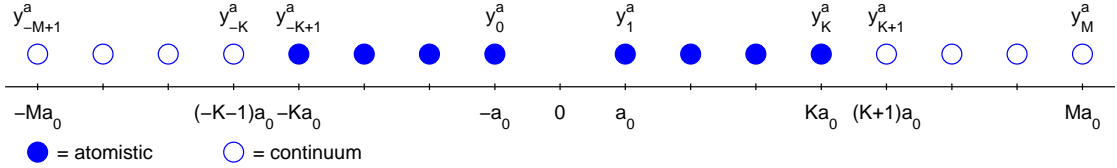}
\caption{Numbering of the atoms. The dislocation is situated in the middle of
  the chain between atoms $y^a_0$ and $y^a_1$.}
\label{FigFKNumbering}
\end{figure}

We study a single defect in the middle of the chain of $2M$ atoms.
To achieve a symmetric description in terms of bonds, we
number the atoms from $-M+1$ to $M$. The defect is situated between the
atoms numbered 0 and 1 (Figure~\ref{FigFKNumbering}).

Recall that the atomistic positions are denoted by $\vec y^a= (y^a_{-M+1},\dots,
y^a_M)\in\R^{2M}$.  The total
potential energy for this atomistic system is then a function $\E^a: \R^{2M} \to
\R$ of the atomistic positions.  For the Frenkel-Kontorova model, the energy,
$\E^a=\E^{a,e}+\E^{a,m}$,
consists of two parts, namely the part which models the elastic energy of the
defect, $\E^{a,e}$, and the part which models the misfit energy on the slip plane, $\E^{a,m}$.

\begin{figure}
\includegraphics[width=0.8\textwidth]{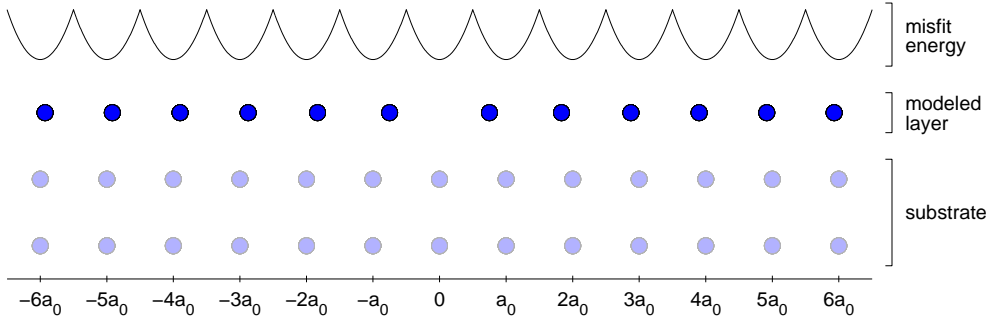}
\caption{Frenkel-Kontorova model. The wells depict the misfit
energy \eqref{EqSubstratePart}.}
\label{FigFKModel}
\end{figure}

The elastic energy is modeled by Hookean (harmonic) springs between nearest-neighbors (NN)
and next-nearest neighbors (NNN), and the total elastic energy is given by
\begin{equation}\label{EqEsFKGeneral}
\E^{a,e}(\vec y^a)=\sum_{i=-M+1}^{M-1} \half k_1 (y^a_{i+1} - y^a_i - a_0)^2
   + \sum_{i=-M+2}^{M-1} \half k_2 (y^a_{i+1} - y^a_{i-1} - 2a_0)^2,
\end{equation}
where the moduli $k_1>0$ and $k_2>0$ describe the
strength of the elastic interactions,
and where $a_0\in\R$ denotes the equilibrium distance.

We note that the asymptotic expansion to second order of any nonlinear
NN/NNN potential energy
\begin{align}\label{true}
 \E(\vec y^a)=
 \sum_{i=-M+1}^{M-1} \psi(y^a_{i+1} - y^a_i)
+ \sum_{i=-M+2}^{M-1} \psi(y^a_{i+1} - y^a_{i-1})
\end{align}
about $ \vec {a}^a=\left[(-M+1) a_0, (-M+2)a_0 , \cdots ,(M-1) a_0 , M a_0\right]^T\in\R^{2M}$
has the form
\begin{equation}\label{bou}
\begin{split}
\E(\vec y^a) & \approx \E(\vec {a}^a)
+\left(\psi'(a_0)+2\psi'(2a_0)\right)\sum_{i=-M+1}^{M-1}
(y^a_{i+1} - y^a_i - a_0)\\
& \qquad -\psi'(2a_0)(y^a_{M} - y^a_{M-1} - a_0)
-\psi'(2a_0)(y^a_{-M+2} - y^a_{-M+1} - a_0)\\
& \qquad +\half \psi''(a_0) \sum_{i=-M+1}^{M-1} (y^a_{i+1} - y^a_i - a_0)^2
   + \half \psi''(2a_0) \sum_{i=-M+2}^{M-1} (y^a_{i+1} - y^a_{i-1} - 2a_0)^2.
\end{split}
\end{equation}
We thus see that the elastic energy \eqref{EqEsFKGeneral}
with $k_1=\psi''(a_0)$ and $k_2=\psi''(2a_0)$ approximates the energy
\eqref{true} to second order
if $\psi'(a_0)+2\psi'(2a_0)=0$ and if we ignore the boundary terms
in the second line of \eqref{bou}.

The misfit energy of the slip plane is modeled
by a periodic potential (Figure~\ref{FigFKModel}). We model this misfit energy by
\begin{equation} \label{EqSubstratePart}
\E^{a,m}(\vec y^a)=\sum_{i=-M+1}^M \half k_0 \left( y^a_i - a_0 \left\lfloor
   \frac{y^a_i}{a_0} + \frac{1}{2} \right\rfloor \right)^2,
\end{equation}
where $\lfloor x \rfloor$ denotes the largest integer smaller than or equal to
$x$, and where the constant $k_0>0$ determines the strength
of the misfit energy.

Altogether, the total potential energy of the atomistic system is given by
\begin{equation}\label{FK}
\begin{split}
\E^a(\vec y^a) & = \E^{a,e}(\vec y^a)+\E^{a,m}(\vec y^a)\\&=
   \half k_1 \sum_{i=-M+1}^{M-1} (y^a_{i+1} - y^a_i - a_0)^2
   + \half k_2 \sum_{i=-M+2}^{M-1} (y^a_{i+1} - y^a_{i-1} - 2a_0)^2\\
&\qquad\qquad+\half k_0 \sum_{i=-M+1}^M \left( y^a_i - a_0 \left\lfloor
   \frac{y^a_i}{a_0} + \frac{1}{2} \right\rfloor \right)^2.
\end{split}
\end{equation}
We restrict ourselves to configurations in which the $M$ leftmost atoms $y_i^a$
for $-M+1 \le i \le 0$ are situated in the interval $\left( \left( i-\frac{3}{2}
  \right) a_0, \left( i-\frac{1}{2} \right) a_0 \right)$, whereas the $M$
rightmost atoms $y_i^a$ for $1 \le i \le M$ are situated in the interval $\left(
  \left( i-\frac{1}{2} \right) a_0, \left( i+\frac{1}{2} \right) a_0 \right)$.
The defect is situated between atoms $y_0^a$ and $y_1^a.$ In this case, the
total energy simplifies to
\begin{equation}\label{EqEnFK}
\begin{split}
  \E^a(\vec y^a) & =
\half k_1\sum_{i=-M+1}^{M-1} (y^a_{i+1} - y^a_i - a_0)^2
  + \half k_2\sum_{i=-M+2}^{M-1} (y^a_{i+1} - y^a_{i-1} - 2a_0)^2
   \\
  & \qquad
+\half k_0 \sum_{i=-M+1}^{0}  \left( y^a_i - (i-1) a_0 \right)^2
  + \half k_0\sum_{i=1}^{M}     \left( y^a_i - i     a_0 \right)^2 .
\end{split}
\end{equation}

%%%%%%%%%%%%%%%%%%%%%%%%%%%%%%%%%%%%%%%%%%%%%%%%%%%%%%%%%%%%%%%%%%%%%%%%%%%%%%%%
\subsection{Quasicontinuum Approximation of the Frenkel-Kontorova Model}

We now apply the quasicontinuum method to the dislocation model described in
Section~\ref{SecFKAtomistic}.  The total energy \eqref{EqEnFK} is split up into
atom-wise contributions, separately for the elastic interactions and the misfit interactions:
\begin{equation} \label{EqDefEaSi}
\begin{split}
  \E^{a,e}_i(\vec y^a) & =
              \fourth k_1(y^a_i    - y^a_{i-1} - a_0)^2
            + \fourth k_1(y^a_{i+1} - y^a_i  - a_0)^2     \\
    & \quad + \fourth k_2(y^a_i    - y^a_{i-2} - 2a_0)^2
            + \fourth k_2(y^a_{i+2} - y^a_i - 2a_0)^2,\\
  \E^{a,m}_i(\vec y^a) & = \begin{cases}
    \half k_0 \left( y^a_i - (i-1) a_0 \right)^2,  & i=-M+1,\ldots,0, \\
    \half k_0 \left( y^a_i - i     a_0 \right)^2,  & i=1,\ldots,M.
  \end{cases}
\end{split}
\end{equation}
To simplify notation, we use the convention that the undefined terms at the
endpoints of the chain are neglected.  We thus have that
\begin{align}\label{to}
\E^a(\vec y^a) = \E^{a,m}(\vec y^a) +\E^{a,e}(\vec y^a) =\sum_{i=-M+1}^M
\left[ \E^{a,m}_i(\vec y^a) + \E^{a,e}_i(\vec y^a) \right].
\end{align}

Since the largest displacement of the atoms is to be expected near the
defect, we deem the atoms $-K+1,\ldots,K$ atomistic and the remaining atoms
$-M+1,\ldots,-K$ and $K+1,\ldots,M$ continuum. Here $K<M$ is some constant whose
optimal value will be determined by the algorithm given in Section~\ref{SecLocNonloc}.

The optimal choice of the repatoms for coarsening is
investigated in the second paper of this series, so we work with a general
formulation which holds for any values of $\ell_j$ for now.  However, there are
two restrictions on the coarsening. Since the atomistic region must not be
coarsened and since we need full refinement in the vicinity of two atoms around
the atomistic region due to the NNN interactions, we have that
\begin{align}
  \ell_j = j, \qquad j=-K-1,\ldots,K+2.
\end{align}
Second, we require that
\begin{equation}\label{require}
  \ell_{-N+1}=-M+1, \qquad
  \ell_{-N+2}=-M+2, \qquad
  \ell_{N-1}=M-1, \qquad \text{and} \qquad
  \ell_{N}=M
\end{equation}
to incorporate the boundary conditions later.

The elastic part $\E^{a,e}_i$ is translationally invariant, so we perform its QC
approximation as described in the previous section. This leads to the continuum
energy density
\begin{align}
      \phi^e(r) & = \half k_1 (r-a_0)^2 + \half k_2(2r-2a_0)^2 \nonumber\\
                & = \half k_{12} (r-a_0)^2
\end{align}
where $k_{12} := k_1 + 4 k_2.$

Regarding the misfit part $\E^{a,m}_i$, the above technique cannot be applied
since the potential is not translationally invariant.  However, there is a
different summation technique to achieve a computationally efficient formulation
which avoids the costly interpolation operator.

To shorten the notation, we let $\sumprime$ indicate the sum in which the first
term and the last term are only counted half:
\begin{align}
  \sumprime_{i=m}^n x_i := \half x_m + \sum_{i=m+1}^{n-1} x_i + \half x_n
\end{align}
where $m<n$ and $x_i\in\R$.  It is easy to verify that
\begin{align} \label{EqPrimeSum}
  \sumprime_{i=0}^m i^2    = \frac{2m^3+m}{6} \qquad \text{and} \qquad
  \sumprime_{i=0}^m i(m-i) = \frac{m^3-m}{6}
\end{align}
for $m>0$.

For all pairs $(j,j+1)$ of continuum repatoms, we now reformulate all terms from
\eqref{EqEcDef} which involve the interaction between $\ell_j$ and $\ell_{j+1}$.
For $j>0$, we get by definition \eqref{EqDefL} of the operator $L$, by
definition \eqref{EqDefEaSi} of $\E^{a,m}_i$, and by \eqref{EqPrimeSum} that
\begin{equation*}
\begin{split}
    \half \sum_{i=\ell_j+1}^{\ell_{j+1}} & \E^{a,m}_i(L^{i-1,i} I \vec y^{qc})
  + \half \sum_{i=\ell_j}^{\ell_{j+1}-1}   \E^{a,m}_i(L^{i,i+1} I \vec y^{qc}) \nonumber\\
& = \half \sum_{i=\ell_j+1}^{\ell_{j+1}}   \E^{a,m}_i(         I \vec y^{qc})
  + \half \sum_{i=\ell_j}^{\ell_{j+1}-1}   \E^{a,m}_i(         I \vec y^{qc}) \nonumber\\
& = \sumprime_{i=\ell_j}^{\ell_{j+1}} \half k_0
   \left( \frac{\ell_{j+1}-i}{\nu_j} y^{qc}_j + \frac{i-\ell_j}{\nu_j} y^{qc}_{j+1} - ia_0 \right)^2
         \nonumber\\
& = \sumprime_{i=\ell_j}^{\ell_{j+1}} \half k_0
   \left( \frac{\ell_{j+1}-i}{\nu_j}(y^{qc}_j-\ell_ja_0)
        + \frac{i-\ell_j}{\nu_j} (y^{qc}_{j+1}-\ell_{j+1}a_0) \right)^2 \nonumber\\
& = \half k_0 \frac{(y^{qc}_j-\ell_ja_0)^2}{\nu_j^2}
    \sumprime_{i=\ell_j}^{\ell_{j+1}} (\ell_{j+1}-i)^2
  + \half k_0 \frac{(y^{qc}_{j+1}-\ell_{j+1}a_0)^2}{\nu_j^2}
    \sumprime_{i=\ell_j}^{\ell_{j+1}} (i-\ell_j)^2 \nonumber\\
& \qquad + k_0 \frac{(y^{qc}_j-\ell_ja_0)(y^{qc}_{j+1}-\ell_{j+1}a_0)}{\nu_j^2}
    \sumprime_{i=\ell_j}^{\ell_{j+1}} (\ell_{j+1}-i)(i-\ell_j) \nonumber\\
& =  \half k_0 \frac{(y^{qc}_j-\ell_ja_0)^2}{\nu_j^2} \frac{2\nu_j^3+\nu_j}{6}
   + \half k_0 \frac{(y^{qc}_{j+1}-\ell_{j+1}a_0)^2}{\nu_j^2} \frac{2\nu_j^3+\nu_j}{6} \nonumber\\
& \qquad + k_0 \frac{(y^{qc}_j-\ell_ja_0)(y^{qc}_{j+1}-\ell_{j+1}a_0)}{\nu_j^2} \frac{\nu_j^3-\nu_j}{6} \nonumber\\
& =: \phi^{m}_{j+}(y^{qc}_j,y^{qc}_{j+1}).
\end{split}
\end{equation*}
For $j<0$, we similarly obtain
\begin{equation*}
\begin{split}
      \half \sum_{i=\ell_j+1}^{\ell_{j+1}} & \E^{a,m}_i(L^{i-1,i} I \vec y^{qc})
    + \half \sum_{i=\ell_j}^{\ell_{j+1}-1}   \E^{a,m}_i(L^{i,i+1} I \vec y^{qc}) \nonumber\\
  & = \half k_0 \frac{(y^{qc}_j-(\ell_j-1)a_0)^2}{\nu_j^2} \frac{2\nu_j^3+\nu_j}{6}
    + \half k_0 \frac{(y^{qc}_{j+1}-(\ell_{j+1}-1)a_0)^2}{\nu_j^2} \frac{2\nu_j^3+\nu_j}{6} \nonumber\\
  & \qquad + k_0 \frac{(y^{qc}_j-(\ell_j-1)a_0)(y^{qc}_{j+1}-(\ell_{j+1}-1)a_0)}{\nu_j^2} \frac{\nu_j^3-\nu_j}{6}
           \nonumber\\
  & =: \phi^{m}_{j-}(y^{qc}_j,y^{qc}_{j+1}).
\end{split}
\end{equation*}

% FOR WHAT DID WE NEED THE FOLLOWING STUFF?
%Note that these formulas are valid within both the atomistic and the continuum
%regions, since the substrate energy is not affected by the continuum approximation.
%For later reference, it is useful to give a formula of this type for the
%substrate energy across the dislocation.  Due to the assumptions $\ell_0=0$ and
%$\ell_1=1,$ we have $\nu_0=1$, so we get for $j=0$:
%\begin{align}
%      \half \sum_{i=\ell_0+1}^{\ell_{1}} & \E^{a,S}_i(L^{i-1,i} I \vec y^{qc})
%    + \half \sum_{i=\ell_0}^{\ell_{1}-1}   \E^{a,S}_i(L^{i,i+1} I \vec y^{qc}) \nonumber\\
%  & = \half k_0 \frac{(y^{qc}_0-(\ell_0-1)a_0)^2}{\nu_0^2} \frac{2\nu_0^3+\nu_0}{6}
%    + \half k_0 \frac{(y^{qc}_1- \ell_1   a_0)^2}{\nu_0^2} \frac{2\nu_0^3+\nu_0}{6} \nonumber\\
%  & =: \phi^{S}_{0}(y^{qc}_0,y^{qc}_1).
%\end{align}

Since $\E^{a,m}_{-K}(L^{-K,-K+1} I \vec y^{qc}) = \E^{a,m}_{-K}(I \vec y^{qc})$
and $\E^{a,m}_{K+1}(L^{K,K+1} I \vec y^{qc}) = \E^{a,m}_{K+1}(I \vec y^{qc})$,
the QC approximation of the chain can be given by
\begin{equation}\label{EqFK_QC}
\begin{split}
  \E^{qc}(\vec y^{qc})
  = &   \sum_{j=-N+1}^{-K} \left[ w_j \phi^e \left( \frac{y^{qc}_{j+1}-y^{qc}_j}{\nu_j} \right)
                                  + \phi^{m}_{j-}(y^{qc}_j,y^{qc}_{j+1}) \right] \\
    & + \half \E^{a,m}_{-K}(I \vec y^{qc})
      + \sum_{j=-K+1}^K \E^a_i(I \vec y^{qc})
      + \half \E^{a,m}_{K+1}(I \vec y^{qc})   \\
    & + \sum_{j=K}^{N-1}   \left[ w_j \phi^e \left( \frac{y^{qc}_{j+1}-y^{qc}_j}{\nu_j} \right)
                                  + \phi^{m}_{j+}(y^{qc}_j,y^{qc}_{j+1}) \right].
 \end{split}
\end{equation}
Note that the interpolation $I\vec y^{qc}$ does not have to be computed
here since the relevant terms only depend on uncoarsened parts of the chain.

Additionally, we will consider the atomistic-continuum approximation
\begin{align} \label{EqFK_C}
\E^{ac}: \R^{2M} \to \R
\end{align}
of the atomistic energy without coarsening. It is given exactly like the QC
approximation \eqref{EqFK_QC} with the only difference being that $\nu_j=1$ and
$\ell_j=j$ everywhere.

%% file: problem.tex
% $Id: problem.tex,v 1.34 2007/04/11 19:50:06 arndt Exp $

%%%%%%%%%%%%%%%%%%%%%%%%%%%%%%%%%%%%%%%%%%%%%%%%%%%%%%%%%%%%%%%%%%%%%%%%%%%%%%%%
%%%%%%%%%%%%%%%%%%%%%%%%%%%%%%%%%%%%%%%%%%%%%%%%%%%%%%%%%%%%%%%%%%%%%%%%%%%%%%%%
%%%%%%%%%%%%%%%%%%%%%%%%%%%%%%%%%%%%%%%%%%%%%%%%%%%%%%%%%%%%%%%%%%%%%%%%%%%%%%%%
\section{Primal and Dual Problems} \label{SecProblem}

%%%%%%%%%%%%%%%%%%%%%%%%%%%%%%%%%%%%%%%%%%%%%%%%%%%%%%%%%%%%%%%%%%%%%%%%%%%%%%%%
\subsection{Problem Setup}

We are now ready to set up the problems we will solve.  We are
interested in finding the minimum of the energy \eqref{EqFK_QC} subject to given
boundary conditions.

We give the boundary conditions by
constraining the deformation of two atoms at
each end of the chain.  This guarantees that the potential with
next-nearest-neighbor interactions can be directly
 applied to all non-boundary atoms
without having to neglect interactions.  We define the spaces
\begin{align}
  V^a := \R^{2M}, \qquad
  V^a_0 := \R^{2M-4}, \qquad
  V^{qc} := \R^{2N}, \qquad
  V^{qc}_0 := \R^{2N-4}.
\end{align}
The spaces $V^a$ and $V^a_0$ will also be used for 
the uncoarsened atomistic-continuum potential $\E^{ac}$,
so there is no need to define spaces $V^{ac}$ and $V^{ac}_0$.
We let ${\vec y}^{bc} \in V^a$  denote any
vector which has the desired boundary values
$y^{bc}_{-M+1},$ $y^{bc}_{-M+2},$ $y^{bc}_{M-1},$ and $y^{bc}_{M},$
and we let ${\vec y}^{bcq} \in V^{qc}$ by any vector satisfying (recall
\eqref{require})
\begin{align} \label{EqBoundaryCond}
y^{bcq}_{-N+1} & = y^{bc}_{-M+1}, &
y^{bcq}_{N-1} & = y^{bc}_{M-1}, \nonumber\\
y^{bcq}_{-N+2} & = y^{bc}_{-M+2}, &
y^{bcq}_{N} & = y^{bc}_{M}.
\end{align}

For any vector $\vec y\in V^a_0$, we denote the
extension by zero boundary conditions to be $J \vec y\in V^a,$ so
\begin{equation}
  J \vec y := \begin{bmatrix} 0 & 0 & \vec y^T & 0 & 0 \end{bmatrix}^T\in\R^{2M},
\end{equation}
and similarly we denote the extension by zero boundary conditions of $\vec y\in V^{qc}_0$
to be $J^{qc} \vec y\in V^{qc}$.  The spaces of admissible solutions are then given
by $ J V^a_0 +{\vec y}^{bc}\subset V^a$ and $  J^{qc}
V^{qc}_0 +{\vec y}^{bcq}\subset V^{qc}$, respectively.  We note that
$J^T:V^a\to V^a_0$ is the restriction operator defined by
\[
(J^T{\vec y})_j=y_j\quad \text{for }j=-M+3,\dots,M-2.
\]

The minima $\bar{\vec y}^a$, $\bar{\vec y}^{ac},$ and $\bar{\vec y}^{qc}$ of
the energy functionals $\E^a$, $\E^{ac},$ and $\E^{qc}$ given by \eqref{EqEnFK},
\eqref{EqFK_C}, and \eqref{EqFK_QC} subject to the above ``clamped'' boundary
conditions are characterized as
\begin{align}
  \bar{\vec y}^a    & := \argmin_{\vec y\in {J V^a_0+\vec y}^{bc}} \E^a(\vec y)
    \in V^a, \\
  \bar{\vec y}^{ac}    & := \argmin_{\vec y\in {J V^a_0+\vec y}^{bc}} \E^{ac}(\vec y)
    \in V^a,  \\
  \bar{\vec y}^{qc} & := \argmin_{\vec y\in {J^{qc} V^{qc}_0+\vec y}^{bcq}} \E^{qc}(\vec y)
    \in V^{qc}.
\end{align}
We note that the minima are uniquely determined because $\E^a$, $\E^{ac}$, and
$\E^{qc}$ are strictly convex.

%%%%%%%%%%%%%%%%%%%%%%%%%%%%%%%%%%%%%%%%%%%%%%%%%%%%%%%%%%%%%%%%%%%%%%%%%%%%%%%%
\subsection{Matrix Formulation} \label{SecMatrixFormulation}

For the subsequent discussion, it will be convenient to reformulate the total
energies in matrix notation:
\begin{subequations}
\begin{align}
  \E^a(\vec y) & = \half (\vec y - \vec a^a)^T D^{aT} E^a D^a (\vec y - \vec a^a)
                 + \half (\vec y - \vec b^a)^T K^a           (\vec y - \vec b^a), \\
  \E^{ac}(\vec y) & = \half (\vec y - \vec a^a)^T D^{aT} E^{ac} D^a (\vec y - \vec a^a)
                 + \half (\vec y - \vec b^a)^T K^a           (\vec y - \vec b^a), \\
  \E^{qc}(\vec y) & = \half (\vec y - \vec a^{qc})^T D^{qcT} E^{qc} D^{qc} (\vec y - \vec a^{qc})
                 + \half (\vec y - \vec b^{qc})^T K^{qc}             (\vec y - \vec b^{qc}).
\end{align}
\end{subequations}
The matrices $D^a \in \R^{(2M-1) \times 2M}$ and $D^{qc} \in \R^{(2N-1)
  \times 2N}$ compute the distance between two adjacent atomistic positions; the
matrices $E^a \in \R^{(2M-1) \times (2M-1)}$, $E^{ac} \in \R^{(2M-1) \times
  (2M-1)}$, and $E^{qc} \in \R^{(2N-1) \times (2N-1)}$ contain
  the spring constants $k_1$, $k_2$, and $k_{12};$
  and the matrices $K^a \in \R^{(2M-1) \times
  (2M-1)}$ and $K^{qc} \in \R^{(2N-1) \times (2N-1)}$ contain
  the misfit constant $k_0$. The vectors $\vec a^a, \vec b^a \in
\R^{2M}$ and $\vec a^{qc}, \vec b^{qc} \in \R^{2N}$ are constants describing
the minimum energy deformations for the elastic energy and misfit energy.  The precise
and lengthy definitions for all of these matrices and vectors
are given in Appendix~\ref{AppMatrixDef}.

If we decompose $\bar{\vec y}^a = J\vec y^a + {\vec y}^{bc}$ for $\vec y^a\in V_0^a,$
then the minimization problem is given as
\begin{align}
  \vec y^a & = \argmin_{\vec y\in V_0^a}  \E^a(J\vec y+{\vec y}^{bc}) \nonumber\\
  & = \argmin_{\vec y\in V_0^a}
     \Big[ \half \big( J\vec y + {\vec y}^{bc} - \vec a^a \big)^T D^{aT} E^a D^a
            \big( J\vec y + {\vec y}^{bc} - \vec a^a \big) \nonumber\\
  & \qquad\qquad\qquad
          + \half \big( J\vec y + {\vec y}^{bc} - \vec b^a \big)^T K^a
            \big( J\vec y + {\vec y}^{bc} - \vec b^a \big)\Big].
\end{align}
We also decompose
$\bar{\vec y}^{ac}
= J\vec y^{ac} + {\vec y}^{bc}$ and $\bar{\vec y}^{qc} = J^{qc}\vec y^{qc} +
{\vec y}^{bcq}$ for $\vec y^{ac}\in V_0^a$ and $\vec y^{qc}\in
V_0^{qc}$, and we then formulate similar minimization problems for
$\vec y^{ac}$ and $\vec y^{qc}.$
Therefore, $\vec
y^a$, $\vec y^{ac},$ and $\vec y^{qc}$ are determined by the linear systems
\begin{subequations} \label{EqPrimal}
\begin{align}
  M^a    \vec y^a   & =  \vec f^a,   \label{EqPrimalA} \\
  M^{ac} \vec y^{ac} & =  \vec f^{ac}, \label{EqPrimalB} \\
  M^{qc} \vec y^{qc} & = \vec f^{qc}, \label{EqPrimalC}
\end{align}
\end{subequations}
where
\begin{equation}
\begin{split} \label{EqPrimal2}
  M^a    & := J^T ( D^{aT} E^a D^{a} + K^a ) J,              \\
  M^{ac} & := J^T ( D^{aT} E^{ac} D^{a} + K^a ) J,            \\
  M^{qc} & := J^{qcT} ( D^{qcT} E^{qc} D^{qc} + K^{qc} ) J^{qc},  \\
  \vec f^a    & := - J^T D^{aT} E^a D^a ({\vec y}^{bc} - \vec a^a)
                   - J^T K^a ({\vec y}^{bc} - \vec b^a), \\
  \vec f^{ac} & := - J^T D^{aT} E^{ac} D^a ({\vec y}^{bc} - \vec a^a)
                   - J^T K^a ({\vec y}^{bc} - \vec b^a), \\
  \vec f^{qc} & := - J^{qcT} D^{qcT} E^{qc} D^{qc} ({\vec y}^{bcq} - \vec a^{qc})
                  - J^{qcT} K^{qc} ({\vec y}^{bcq} - \vec b^{qc}).
\end{split}
\end{equation}
We note that the matrices $M^a$, $M^{ac},$ and $M^{qc}$ are positive definite, so
the total energies admit a single global minimum and no other local minimum.

%%%%%%%%%%%%%%%%%%%%%%%%%%%%%%%%%%%%%%%%%%%%%%%%%%%%%%%%%%%%%%%%%%%%%%%%%%%%%%%%
\subsection{Goal-Oriented Error Estimation}

To compare the approximate QC model to the original atomistic model, we have to
analyze how much the solution $\vec y^a$ of the atomistic model deviates from
the solution $\vec y^{qc}$ of the QC model.  This deviation, which can be viewed
as an approximation error, can be measured in different ways, for example as $\|\vec
y^a - J^T I J^{qc} \vec y^{qc}\|$ for some norm $\|\vec \cdot \|$. Here we follow
a different approach, namely we measure the error of a {\em quantity of
  interest} denoted by $Q(\vec y)$ for some function $Q:\R^{2M-4}\to\R$.
Hence, we intend to estimate
\begin{align}
  Q(\vec y^a) - Q(J^T I J^{qc} \vec y^{qc}).
\end{align}
We will assume for simplicity that $Q$ is linear and thus
has a representation $Q(\vec y) = \vec q^T \vec y$ for some vector $\vec
q \in V^a_0$.

For our application, a natural quantity of interest is the size of the
dislocation, that is, the distance between the two atoms $y_0$ and $y_1$ to the left
and right of the dislocation. This gives us
\begin{align}
  Q(\vec y) = \vec q^T \vec y = y_1 - y_0 \qquad \text{with} \qquad
  \vec q = [0,\ldots,0,-1,1,0,\ldots,0]^T.
\end{align}

Two different sources of error arise during the QC approximation, namely the
localization of the potential energy, that is, the passage from the atomistic to the continuum
formulation on the one hand, and the coarsening in the continuum region by the
restriction to the repatoms on the other hand. We denote these two errors by
\begin{align}\label{pri}
  \vec e := \vec y^a - \vec y^{ac}
  \qquad \text{and} \qquad
  \vec e^{acqc} := \vec y^{ac} - J^T I J^{qc} \vec y^{qc}.
\end{align}

It makes sense to study these sources independently.  Employing the linearity of
$Q$, we have that
\begin{align} \label{EqErrorPrimal}
  |Q(\vec y^a) - Q( J^T I J^{qc} \vec y^{qc})|
  =   |Q(\vec e)  +  Q(\vec e^{acqc})|
  \le |Q(\vec e)| + |Q(\vec e^{acqc})|.
\end{align}
The error term $|Q(\vec e)|$ will be studied in
Section~\ref{SecLocNonloc}, and the error term $|Q(\vec e^{acqc})|$ will be studied in
the second part of this paper series.

%%%%%%%%%%%%%%%%%%%%%%%%%%%%%%%%%%%%%%%%%%%%%%%%%%%%%%%%%%%%%%%%%%%%%%%%%%%%%%%%
\subsection{Dual Problems}

To facilitate the goal-oriented error analysis, we introduce the dual problems
\begin{subequations} \label{EqDual}
\begin{align}
  M^a    \vec g^a    & =     \vec q,  \label{EqDualA} \\
  M^{ac} \vec g^{ac} & =     \vec q, \label{EqDualB} \\
  M^{qc} \vec g^{qc}  & = J^{qcT} I^T J \vec q, \label{EqDualC}
\end{align}
\end{subequations}
for $\vec g^a, \vec g^{ac} \in \R^{2M-4},$ and $\vec g^{qc} \in \R^{2N-4}$.  We note
that the dual problems differ from the primal problems only by the right hand
side since the matrices $M^a$, $M^{ac},$ and $M^{qc}$ are symmetric.

The solutions $\vec g^a$, $\vec g^{ac}$ and $\vec g^{qc}$ can be viewed as
influence functions: They describe how the error at a specific point in the
domain influences the error measured in terms of the goal function.

Analogously to the primal errors \eqref{pri}, we define the dual
errors
\begin{align} \label{EqErrorDual}
  \hat{\vec e} := \vec g^a - \vec g^{ac}
  \qquad \text{and} \qquad
  \hat{\vec e}^{acqc} := \vec g^{ac} - J^T I J^{qc} \vec g^{qc}.
\end{align}

In addition, we will need the primal and dual residuals
\begin{align}
  R^a(\vec y) & := M^a \left( \vec y^a - \vec y \right)
                 = \vec f^a - M^a \vec y, \nonumber\\
  R^{ac}(\vec y) & := M^{ac} \left( \vec y^{ac} - \vec y \right)
                 = \vec f^{ac} - M^{ac} \vec y, \nonumber\\
  \hat{R}^a(\vec g) & := M^a \left( \vec g^a - \vec g \right)
      = \vec q - M^a \vec g, \nonumber\\
  \hat{R}^{ac}(\vec g) & := M^{ac} \left( \vec g^{ac} - \vec g \right)
      = \vec q - M^{ac} \vec g.
\end{align}

%% file: errorcont.tex
% $Id: errorcont.tex,v 1.36 2007/04/11 19:50:07 arndt Exp $

%%%%%%%%%%%%%%%%%%%%%%%%%%%%%%%%%%%%%%%%%%%%%%%%%%%%%%%%%%%%%%%%%%%%%%%%%%%%%%%%
%%%%%%%%%%%%%%%%%%%%%%%%%%%%%%%%%%%%%%%%%%%%%%%%%%%%%%%%%%%%%%%%%%%%%%%%%%%%%%%%
%%%%%%%%%%%%%%%%%%%%%%%%%%%%%%%%%%%%%%%%%%%%%%%%%%%%%%%%%%%%%%%%%%%%%%%%%%%%%%%%
\section{Error Estimation for Atomistic vs.\ Continuum Modeling}  \label{SecLocNonloc}

In this section, we estimate the error $|Q(\vec e)|$ arising from the
approximation of an atomistic model by a continuum model.  We consider
$\vec y^{ac}$ and $\vec g^{ac}$ to be computable, although in practice we can
only compute the coarsened approximations $\vec y^{qc}$ and $\vec g^{qc}.$

To this end, we
adapt a technique introduced in \cite{OdenPrudhomme:2002} and
\cite{OdenVemaganti:2000} to estimate the modeling error for an elasticity model
with rapidly oscillating coefficients and its homogenized version. We generalize
this technique such that it allows for different right hand sides $\vec f^a$ and
$\vec f^{ac}$ of the primal problem \eqref{EqPrimal} instead of a common right
hand side as it is used in the above-mentioned works.

We have
\begin{align} \label{EqBasicDual}
  Q(\vec y^a) - Q(\vec y^{ac})
  & = \vec q^T \vec e = \vec g^{aT} M^a \vec e
    = (\vec g^{acT} + \hat{\vec e}) M^a \vec e \nonumber\\
  & = \vec g^{acT} R^a(\vec y^{ac}) + \hat{\vec e}^T M^a \vec e.
\end{align}
The term $\vec g^{acT} R^a(\vec y^{ac})$ can be computed, whereas $\hat{\vec e}^T M^a
\vec e$ cannot because both $\vec e$ and $\hat{\vec e}$ are numerically unknown.
Instead, we estimate $\hat{\vec e}^T M^a \vec e$ from above and from below by
quantities that actually can be computed.

We will give two different error estimators $\eta_1$ and $\eta_2$. Before, we
need to derive some auxiliary estimates to facilitate their development and
analysis.

%%%%%%%%%%%%%%%%%%%%%%%%%%%%%%%%%%%%%%%%%%%%%%%%%%%%%%%%%%%%%%%%%%%%%%%%%%%%%%%
\subsection{Auxiliary Estimates}

We reformulate the difference $\vec y^a - \vec y^{ac}$ of the
respective solutions in terms of a difference of the energy
matrices. To this end, we define the perturbation matrix
\begin{align} \label{EqI0Def}
  P := {\mathcal I} - (E^a)^{-1} E^{ac}
\end{align}
where ${\mathcal I}$ denotes the identity matrix. Note that $E_a P = E^a -
E^{ac}$.

\begin{lemma} \label{LemmaI0}
  For any $\alpha, \beta\in\R,$ we have that
  \begin{align}
    M^a (\alpha\vec e + \beta\hat{\vec e})
    = - J^T D^{aT} E^a P D^a \big[ \alpha (J\vec y^{ac} + \vec y^{bc} - \vec a^a)
      + \beta J \vec g^{ac}  \big].
  \end{align}
\end{lemma}

\begin{proof}
We conclude from \eqref{EqPrimal} that
\begin{align}
M^a \vec e & =M^a\vec y^a-M^{ac}\vec y^{ac}
+\left(M^{ac}-M^a\right)\vec y^{ac} \nonumber\\
&=\vec f^a-\vec f^{ac}
+\left(M^{ac}-M^a\right)\vec y^{ac},
\end{align}
and similarly since $M^a\vec g^a=M^{ac}\vec g^{ac}=\vec q$ that
\begin{align}
M^a \hat{\vec e} = M^a(\vec g^a-\vec g^{ac}) = (M^{ac}-M^a)\vec g^{ac}.
\end{align}
Thus, it follows from \eqref{EqPrimal2} and \eqref{EqI0Def} that
\begin{equation}
\begin{split}
M^a (\alpha\vec e + \beta\hat{\vec e})
%& = J^T D^{aT} E^a D^a J \left[ \alpha (\vec y^a - \vec y^{ac})
%               + \beta (\vec g^a - \vec g^{ac}) \right] \nonumber\\
& = \alpha \left[(M^{ac}-M^a)\vec y^{ac}+\vec f^a-\vec f^{ac}\right]
    + \beta ( M^{ac}-M^a ) \vec g^{ac} \\
& = J^T D^{aT} (E^{ac}-E^a) D^a \big[ \alpha ( J \vec y^{ac} + \vec y^{bc} - \vec a^a )
               + \beta J \vec g^{ac} \big] \\
& = - J^T D^{aT} E^a P D^a \big[ \alpha ( J \vec y^{ac} + \vec y^{bc} - \vec a^a )
               + \beta J \vec g^{ac} \big].
\end{split}
\end{equation}
We note that the $K^a$-related terms cancel here, because they coincide
for the atomistic model and the continuum model.
\end{proof}

\begin{lemma} \label{Lemma2}
  We have that
  \begin{align}
    \| \alpha \vec e + \beta \hat{\vec e} \|_{M^a}
    \le \big\| P D^a \big[ \alpha ( J \vec y^{ac} + \vec y^{bc} - \vec a^a )
       + \beta J \vec g^{ac} \big] \big\|_{E^a}.
  \end{align}
\end{lemma}

We note that the right hand side is numerically computable.

\begin{proof}
  To shorten the notation, we abbreviate $\vec z = \alpha ( J \vec y^{ac} + \vec y^{bc} - \vec
  a^a ) + \beta J \vec g^{ac}$. By Lemma~\ref{LemmaI0}, we have
  \begin{equation}
\begin{split}
    \| \alpha \vec e + \beta \hat{\vec e} \|_{M^a}
    & = \sup_{\vec v \in V^a_0 \setminus \{0\}}
        \frac{\vec v^T M^a (\alpha \vec e + \beta \hat{\vec e})}{\|\vec v\|_{M^a}} \\
    & = \sup_{\vec v \in V^a_0 \setminus \{0\}}
        \frac{- \vec v^T J^T D^{aT} E^a P D^a \vec z}{\|\vec v\|_{M^a}} \\
    & \le \sup_{\vec v \in V^a_0 \setminus \{0\}}
        \frac{\| D^a J \vec v \|_{E^a} \| P D^a \vec z\|_{E^a}}{\|D^a J \vec v\|_{E^a}}\\
    & = \| P D^a \vec z\|_{E^a}.
\end{split}
\end{equation}
  Here we have used that $\|D^a J \vec v\|_{E^a} \le \|\vec v\|_{M^a}$ because the matrix
  $K^a$ in \eqref{EqPrimal2} is positive definite.
\end{proof}

%%%%%%%%%%%%%%%%%%%%%%%%%%%%%%%%%%%%%%%%%%%%%%%%%%%%%%%%%%%%%%%%%%%%%%%%%%%%%%%
\subsection{First Error Estimator}

We are now ready to derive the first error estimator, $\eta_1$.
By the parallelogram identity, we have for all $\sigma \ne 0$ that
\begin{equation}\label{EqParallelogram}
\begin{split}
  \hat{\vec e}^T M^a \vec e
  & = (\sigma^{-1}\hat{\vec e}^T) M^a (\sigma\vec e) \\
  & = \textstyle \frac{1}{4} \| \sigma\vec e + \sigma^{-1}\hat{\vec e} \|_{M^a}^2
               - \frac{1}{4} \| \sigma\vec e - \sigma^{-1}\hat{\vec e} \|_{M^a}^2.
\end{split}
\end{equation}
In the following, we will determine computable constants $\eta_\text{low}^+$, $\eta_\text{low}^-$,
$\eta_\text{upp}^+$ and $\eta_\text{upp}^-$ such that
 \begin{equation}\label{eta}
\begin{split}
  \eta_\text{low}^+ \le \| \sigma\vec e + \sigma^{-1}\hat{\vec e} \|_{M^a} \le \eta_\text{upp}^+,
  \\
  \eta_\text{low}^- \le \| \sigma\vec e - \sigma^{-1}\hat{\vec e} \|_{M^a} \le
  \eta_\text{upp}^-.
\end{split}
\end{equation}
From Lemma~\ref{Lemma2}, we immediately get the upper estimates
$\eta_\text{upp}^+$ and $\eta_\text{upp}^-$:
 \begin{equation}\label{EqEtaUppDef}
\begin{split}
  \eta_\text{upp}^+
    & := \big\| P D^a \big[ \sigma ( J \vec y^{ac} + \vec y^{bc} - \vec a^a )
       + \sigma^{-1} J \vec g^{ac} \big] \big\|_{E^a}, \\
   \eta_\text{upp}^-
    & := \big\| P D^a \big[ \sigma ( J \vec y^{ac} + \vec y^{bc} - \vec a^a )
       - \sigma^{-1} J \vec g^{ac} \big] \big\|_{E^a}.
\end{split}
\end{equation}
We note that $\eta_\text{low}^+$, $\eta_\text{low}^-$,
$\eta_\text{upp}^+$ and $\eta_\text{upp}^-$ will depend on $\sigma$, but
the estimates \eqref{eta} will hold for any $\sigma \ne 0$.  We will now choose $\sigma$ in such a
way that the estimates are as sharp as possible, that is, such that $\eta_\text{upp}^+$ and
$\eta_\text{upp}^-$ are smallest.

\begin{lemma}
  Both $\eta_\text{upp}^+$ and $\eta_\text{upp}^-$ given by \eqref{EqEtaUppDef}
  attain their minima for
  \begin{align} \label{EqS}
    \bar\sigma := \sqrt{\frac{\| P D^a J \vec g^{ac} \|_{E^a}}
                        {\| P D^a ( J \vec y^{ac} + \vec y^{bc} - \vec a^a ) \|_{E^a}}}.
  \end{align}
\end{lemma}

\begin{proof}
  We have that
  \begin{align}
    (\eta_\text{upp}^\pm)^2
    & = \sigma^2 \big\| P D^a ( J \vec y^{ac} + \vec y^{bc} - \vec a^a ) \big\|_{E^a}^2
    \pm 2 \vec g^{acT} J^T D^{aT} P^T E^a P D^a ( J \vec y^{ac} + \vec y^{bc} - \vec a^a )
    \nonumber\\
    & \qquad + \sigma^{-2} \big\| P D^a J \vec g^{ac} \big\|_{E^a}^2.
  \end{align}
  Setting the first derivative of the mapping $\sigma \mapsto  (\eta_\text{upp}^\pm)^2$ to zero,
  we obtain the condition
  \begin{align}
    2\bar\sigma \big\| P D^a ( J \vec y^{ac} + \vec y^{bc} - \vec a^a ) \big\|_{E^a}^2
    - 2\bar\sigma^{-3} \big\| P D^a J \vec g^{ac} \big\|_{E^a}^2 = 0
  \end{align}
  for critical points of $(\eta_\text{upp}^\pm)^2$.  This equation has the
  unique positive solution \eqref{EqS}. Because $\lim_{|\sigma|\to\infty} \eta_\text{upp}^\pm
  = \lim_{\sigma\to0} \eta_\text{upp}^\pm = \infty$, this point corresponds to a
  minimum. Hence the quantities $\eta_\text{upp}^\pm$ attain their minima at
  $\sigma=\bar\sigma$.
\end{proof}

Regarding the lower bounds $\eta_\text{low}^+$ and $\eta_\text{low}^-$, we have
\begin{equation}
\begin{split}
  \| \bar\sigma\vec e \pm \bar\sigma^{-1}\hat{\vec e} \|_{M^a}
  & = \sup_{\vec v \in V^a_0 \setminus \{0\}}
      \frac{\vec v^T M^a (\bar\sigma\vec e \pm \bar\sigma^{-1}\hat{\vec e})}{\|\vec v\|_{M^a}} \\
  & = \sup_{\vec v \in V^a_0 \setminus \{0\}}
      \frac{\vec v^T (\bar\sigma R^a(\vec y^{ac}) \pm \bar\sigma^{-1} \hat R^a(\vec g^{ac}))}{\|\vec v\|_{M^a}}\\
  & \ge \frac{\vec v_0^T (\bar\sigma R^a(\vec y^{ac}) \pm \bar\sigma^{-1} \hat R^a(\vec g^{ac}))}{\|\vec v_0\|_{M^a}}
\end{split}
\end{equation}
for any vector $\vec v_0 \in V^a_0 \setminus \{0\}$. Numerically, we have the
two vectors $\vec y^{ac}$ and $\vec g^{ac}$ at our disposal, hence it makes sense to
take a linear combination $\vec v_0 = \vec y^{ac} + \theta^\pm \vec g^{ac}$.  Here we
follow the strategy of \cite{OdenVemaganti:2000} and choose $\theta^\pm$ as the
critical points of $\eta_\text{low}^\pm$.

\begin{lemma}
  Let
  \begin{align}
    \vec r^\pm = \bar\sigma R^a(\vec y^{ac}) \pm \bar\sigma^{-1} \hat R^a(\vec g^{ac}).
  \end{align}
  Then the lower bounds
  \begin{align}
    \eta_\text{low}^\pm
    :=\frac{(\vec y^{ac} + \theta^\pm \vec g^{ac})^T \vec r^\pm}
           {\| \vec y^{ac} + \theta^\pm \vec g^{ac} \|_{M^a}}
  \end{align}
  have a unique critical point for
  \begin{align}
    \bar \theta^\pm :=
    \frac{ \vec r^{\pm T} \vec y^{ac} \; \vec g^{acT} M^a \vec y^{ac} -
           \vec r^{\pm T} \vec g^{ac} \; \|\vec y^{ac}\|_{M^a}^2}
         { \vec r^{\pm T} \vec g^{ac} \; \vec g^{acT} M^a \vec y^{ac} -
           \vec r^{\pm T} \vec y^{ac} \; \|\vec g^{ac}\|_{M^a}^2}.
  \end{align}
\end{lemma}

\begin{proof}
  We have
  \begin{align}
    \frac{\D}{\D \theta^\pm} \eta_\text{low}^\pm =
    \frac{ \vec r^{\pm T} \vec g^{ac} \| \vec y^{ac} + \theta^\pm \vec g^{ac} \|_{M^a} -
           \vec r^{\pm T} (\vec y^{ac} + \theta^\pm \vec g^{ac})
           \frac{ (\vec y^{ac} + \theta^\pm \vec g^{ac})^T M^a \vec g^{ac} }
                { \| \vec y^{ac} + \theta^\pm \vec g^{ac} \|_{M^a} } }
         { \| \vec y^{ac} + \theta^\pm \vec g^{ac} \|_{M^a}^2 }.
  \end{align}
  Setting this expression to zero and solving for $\theta^\pm$ leads to the above
  condition.
\end{proof}

However, let us note that this critical point is not necessarily a maximum of
$\eta_\text{low}^\pm,$ which would be optimal for bound \eqref{eta}. Depending on
the actual vectors $\vec y^{ac}$ and $\vec g^{ac}$, it can be shown that this critical
point could be a minimum.

Now we have all necessary ingredients to construct the error estimator
$\eta_1$.  From \eqref{EqBasicDual} and \eqref{EqParallelogram}, we get the
computable estimate
\begin{align} \label{EqEstLowUpp}
\textstyle
\vec g^{acT} R^a(\vec y^{ac}) + \frac{1}{4} (\eta_\text{low}^+)^2
                      - \frac{1}{4} (\eta_\text{upp}^-)^2
\le Q(\vec y^a) - Q(\vec y^{ac}) \le
\vec g^{acT} R^a(\vec y^{ac}) + \frac{1}{4} (\eta_\text{upp}^+)^2
                      - \frac{1}{4} (\eta_\text{low}^-)^2.
\end{align}

At first sight, this looks like we could get an estimate for $|Q(\vec y^a) -
Q(\vec y^{ac})|$ from both above and below.
However, this is only true if both the left hand side
and the right hand side have the same sign, which in general does not hold.  But
we get the following estimate.

\begin{theorem} \label{Theorem1}
  We have that
  \begin{align}
    \left| Q(\vec y^a) - Q(\vec y^{ac}) \right| \le \eta_1,
  \end{align}
  where the computable error estimator is defined as
  \begin{align}  \label{EqEta1Def}
    \textstyle \eta_1 :=
    \max\left( \left| \vec g^{acT} R^a(\vec y^{ac}) + \frac{1}{4} (\eta_\text{low}^+)^2
                                             - \frac{1}{4} (\eta_\text{upp}^-)^2 \right|,
      \left| \vec g^{acT} R^a(\vec y^{ac}) + \frac{1}{4} (\eta_\text{upp}^+)^2
                                             - \frac{1}{4} (\eta_\text{low}^-)^2 \right| \right).
  \end{align}
\end{theorem}

We note that the computation of the $\eta_{upp}^\pm$ terms involves the solution
of a linear system with matrix $E^a$ as its inverse appears in the operator
$P$. The matrix $E^a$ is not diagonal, but has condition number $\Oh(1)$. So
this is negligible compared to what would be necessary to solve the original
atomistic problem which includes the operator $D^{aT} E^a D^a$ with condition
number $\Oh(M^2)$.

%%%%%%%%%%%%%%%%%%%%%%%%%%%%%%%%%%%%%%%%%%%%%%%%%%%%%%%%%%%%%%%%%%%%%%%%%%%%%%%
\subsection{Second Error Estimator}

There is no reasonable way to decompose the error estimator $\eta_1$ into a sum of
element-wise or atom-wise contributions due to the $\eta_{low}^\pm$ terms.
Therefore, we derive another error estimator $\eta_2$ which allows for such a
decomposition, at the price of a less accurate estimate than $\eta_1$.

\begin{theorem}
  We have that
  \begin{align}
    \left| Q(\vec y^a) - Q(\vec y^{ac}) \right| \le \eta_2
    \le \sum_{i=-M+3}^{M-2} \eta_{2,i}^{at} + \sum_{i=-M+1}^{M-1} \eta_{2,i}^{el}
  \end{align}
  where the computable global error estimator $\eta_2$ and the computable local
  error estimators, $\eta_{2,i}^{at}$ and $\eta_{2,i}^{el}$, associated with atoms and
  elements, respectively, are defined as
  \begin{equation} \label{EqEta2Def}
  \begin{tabular}{r@{}ll}
    $\eta_2$     & \multicolumn{2}{@{}l}{$\; := \left| \vec g^{acT} R^a(\vec y^{ac}) \right|
                 + \|P D^a (J\vec y^{ac} + \vec y^{bc} - \vec a^a)\|_{E^a}
                   \|P D^a J\vec g^{ac}\|_{E^a},$} \\[.5em]
    $\eta_{2,i}^{at}$ & $\; := \left| g^{ac}_i R^a(\vec y^{ac})_i \right|,$
                 & $\qquad i=-M+3,\ldots,M-2,$ \\[.5em]
    $\eta_{2,i}^{el}$ & \multicolumn{2}{@{}l}{$
                 \;:= \half \left| \big( P D^a (J\vec y^{ac} + \vec y^{bc} - \vec a^a) \big)_i
                      \big( (E^a-E^{ac}) D^a (J\vec y^{ac} + \vec y^{bc} - \vec a^a) \big)_i 
                      \right|$} \\[.5em]
               & $\qquad + \half \left| (P D^a J\vec g^{ac})_i 
                    \big( (E^a-E^{ac}) D^a J\vec g^{ac} \big)_i \right|,$
               & $\qquad i=-M+1,\ldots,M-1.$
  \end{tabular}
  \end{equation}
\end{theorem}

\begin{proof}
  From \eqref{EqBasicDual} and Lemma~\ref{Lemma2}, we conclude that
 \begin{equation}\label{EqEst2Ineq}
  \begin{split}
    \left| Q(\vec y^a) - Q(\vec y^{ac}) \right|
    & \le \left| \vec g^{acT} R^a(\vec y^{ac}) \right|
        + \left| \hat{\vec e}^T M^a \vec e \right| \\
    & \le \left| \vec g^{acT} R^a(\vec y^{ac}) \right|
        + \| \hat{\vec e} \|_{M^a} \| \vec e \|_{M^a} \\
    & \le \left| \vec g^{acT} R^a(\vec y^{ac}) \right|
        + \|P D^a (J\vec y^{ac} + \vec y^{bc} - \vec a^a)\|_{E^a}
          \|P D^a J\vec g^{ac}\|_{E^a} \\
    & = \eta_2,
  \end{split}
  \end{equation}
  which gives us the global estimate. For the decomposition into local
  contributions, we further estimate
\begin{equation}\label{EqEst2Splitting}
  \begin{split}
    \eta_2 & = \left| \vec g^{acT} R^a(\vec y^{ac}) \right|
             + \|P D^a (J\vec y^{ac} + \vec y^{bc} - \vec a^a)\|_{E^a}
             \|P D^a J\vec g^{ac}\|_{E^a}
              \\
    & \le      \left| \vec g^{acT} R^a(\vec y^{ac}) \right|
           + \half \|P D^a (J\vec y^{ac} + \vec y^{bc} - \vec a^a)\|_{E^a}^2
           + \half \|P D^a J\vec g^{ac}\|_{E^a}^2 \\
   & \le \sum_{i=-M+3}^{M-2} \left| g^{ac}_i R^a(\vec y^{ac})_i \right|
      +  \half \sum_{i=-M+1}^{M-1} \left|
           (P D^a J\vec g^{ac})_i \big( (E^a-E^{ac}) D^a J\vec g^{ac} \big)_i \right| \\
   & \qquad +  \half\sum_{i=-M+1}^{M-1} \left|
            \big(P D^a (J\vec y^{ac} + \vec y^{bc} - \vec a^a) \big)_i
           \big( (E^a-E^{ac}) D^a (J\vec y^{ac} + \vec y^{bc} - \vec a^a) \big)_i \right| \\
  & = \sum_{i=-M+3}^{M-2} \eta_{2,i}^{at} + \sum_{i=-M+1}^{M-1} \eta_{2,i}^{el},
  \end{split}
  \end{equation}
  which completes the proof.
\end{proof}

Let us remark that instead of the first inequality in \eqref{EqEst2Splitting},
one can get an apparently better estimate
\begin{equation}
  \begin{split}
  \lefteqn{\|P D^a (J\vec y^{ac} + \vec y^{bc} - \vec a^a)\|_{E^a} \|P D^a J\vec g^{ac}\|_{E^a}}
   \\
  & \qquad \qquad = \half \gamma  \|P D^a (J\vec y^{ac} + \vec y^{bc} - \vec a^a)\|^2_{E^a}
    + \half \gamma^{-1} \|P D^a J\vec g^{ac}\|^2_{E^a}
\end{split}
  \end{equation}
by introducing the additional weight factor
\begin{align}
  \gamma := \frac{\|P D^a J\vec g^{ac}\|_{E^a}}
              {\|P D^a (J\vec y^{ac} + \vec y^{bc} - \vec a^a)\|_{E^a}},
\end{align}
and then decomposing the resulting terms similar to the above.  However, our
numerical results showed that this modification does not significantly improve
the decomposed error estimator for the application considered here.

%% file: numericscont.tex
% $Id: numericscont.tex,v 1.24 2007/04/12 15:41:16 arndt Exp $

%%%%%%%%%%%%%%%%%%%%%%%%%%%%%%%%%%%%%%%%%%%%%%%%%%%%%%%%%%%%%%%%%%%%%%%%%%%%%%%%
%%%%%%%%%%%%%%%%%%%%%%%%%%%%%%%%%%%%%%%%%%%%%%%%%%%%%%%%%%%%%%%%%%%%%%%%%%%%%%%%
%%%%%%%%%%%%%%%%%%%%%%%%%%%%%%%%%%%%%%%%%%%%%%%%%%%%%%%%%%%%%%%%%%%%%%%%%%%%%%%%
\section{Numerics}  \label{SecNumericalResults}

In the preceding sections, we constructed the error estimators
$\eta_1$ and $\eta_2$. We will now give an algorithm for
adaptive atomistic-continuum modeling based on these error
estimators. Then we will present
and discuss some numerical results.

%%%%%%%%%%%%%%%%%%%%%%%%%%%%%%%%%%%%%%%%%%%%%%%%%%%%%%%%%%%%%%%%%%%%%%%%%%%%%%%%
\subsection{Algorithm} \label{SecAlgorithm}

The error estimator $\eta_1$ should give a better estimate of the error
than $\eta_2$, because $\eta_2$ involves the inequality $\left| \hat{\vec
    e}^T M^a \vec e \right| \le \| \hat{\vec e} \|_{M^a} \| \vec e \|_{M^a}$ in
\eqref{EqEst2Ineq} in contrast to the parallelogram identity for $\eta_1$.
However, $\eta_2$ can be decomposed into atom-wise and element-wise contributions
$\eta_{2,i}^{at}$ and $\eta_{2,i}^{el}$, whereas the $\eta_{low}^\pm$ terms in
$\eta_1$ do not admit a reasonable decomposition that can be used for
atomistic-continuum adaptivity.

We make use of this by employing the sharper estimate $\eta_1$ to determine
whether a given global error tolerance $\tau_{gl}$ for the error in an adaptive
algorithm has already been achieved or not. If not, we use the decomposed
estimates $\eta_{2,i}^{at}$ and $\eta_{2,i}^{el}$ to determine where the more
precise atomistic modeling is needed.  This leads us to the following algorithm:

\begin{enumerate}
\item[(1)] Choose $\tau_{gl}$. Model all atoms as a continuum. Set $\tau_{at} \leftarrow \tau_{gl}$.
\item[(2)] Solve primal problem \eqref{EqPrimalB} for $\vec y^{ac}$ and dual problem
           \eqref{EqDualB} for $\vec g^{ac}$.
\item[(3)] Compute error estimator $\eta_1$ from \eqref{EqEta1Def}.
\item[(4)] If $\eta_1 \le \tau_{gl},$ then stop.
\item[(5)] Compute local error estimators $\eta_{2,i}^{at}$ and $\eta_{2,i}^{el}$
           from \eqref{EqEta2Def}.
\item[(6)] Set $\tau_{at} \leftarrow \frac{\tau_{at}}{\tau_{div}}$.
\item[(7)] Make all atoms $i$ atomistic for which
  \begin{align} \label{EqAdaptCrit}
    \eta_{2,i}^{tot}:=\eta_{2,i}^{at} + \half \big( \eta_{2,i-1}^{el} + \eta_{2,i}^{el} \big) \ge \tau_{at}.
  \end{align}
\item[(8)] Go to (2).
\end{enumerate}
Here $\tau_{div}>1$ is a constant factor which describes how fast the atom-wise
tolerance $\tau_{at}$ should decrease during adaption. Our experience has been that
$\tau_{div}=10$ is a reasonable choice.

The crucial adaption step is (7). The adaption criterion \eqref{EqAdaptCrit}
deems all atoms to be modeled atomistically if the associated error from the
decomposition of $\eta_2$ exceeds the atomistic error tolerance $\tau_{at}$. Here
the element-wise errors $\eta_{2,i}^{el}$ are distributed equally to the two
adjacent atoms $i$ and $i+1$.

For the dislocation at the center of the chain and the chosen goal function, we
expect that the atomistic repatoms always form a symmetric interval around the center.
We have used the above adaptive atomistic-continuum
algorithm to approximate our Frenkel-Kontorova model
and have always found that the atomistic region is the set
of atoms $-K+1, \ldots, K$ for some $K$ depending on $M$ and
$\tau_{gl}.$
Thus, the modeling approach given in Section~\ref{SecFKModel}
of restricting to an atomistic region consisting of atoms $-K+1, \ldots, K$ for some $K$
rather than considering a more general atomistic region is justified {\em a posteriori}.

\begin{comment}
For the special case of our Frenkel-Kontorova model, recall that the atoms
$-K+1, \ldots, K$ are modeled atomistically and that the remaining atoms are
continuum atoms. Therefore, we initially choose $K=0$ to get a fully continuum
chain.  In the adaption step (7), we increase $K$ such that the interval $-K+1,
\ldots, K$ comprises all atoms deemed atomistic by the error estimator:
\begin{align} \label{EqK}
  K \leftarrow \max \big( K, \max\big\{ j\in\{1,\ldots,M\}:
      \eqref{EqAdaptCrit} \text{ holds for } i=j \text{ or for }
                 i=-j+1 \big\} \big).
\end{align}

For the dislocation at the center of the chain and the chosen goal function, we
expect that the atomistic repatoms always form an interval around the center such
that \eqref{EqK} does not tag any repatoms atomistic which need not be atomistic.
The numerical results below confirm this anticipated behavior. Let us note that
tagging atoms unnecessarily to be atomistic would degrade the performance, but
would otherwise be harmless.
\end{comment}
%%%%%%%%%%%%%%%%%%%%%%%%%%%%%%%%%%%%%%%%%%%%%%%%%%%%%%%%%%%%%%%%%%%%%%%%%%%%%%%%
\subsection{Numerical Results} \label{SecNumResults}

The algorithm has been implemented as described above.
The boundary conditions were chosen as
\begin{align}
  y^{bc}_{-M+1} = -M, \qquad
  y^{bc}_{-M+2} = -M+1, \qquad
  y^{bc}_{M-1} = M-1, \qquad
  y^{bc}_{M} = M.
\end{align}
%compare the definition \eqref{EqBoundaryCond}.
The elastic constants are $k_0=1$ and $k_1=k_2=2$.

\begin{table}
\begin{tabular}{|c|c|c|c|c|}
\hline
       M & iteration &   K & $\tau_{at}$ & $\eta_1$ \\
\hline
     100 &   1 &   0 & 1.000000e-10 & 3.899207e-02 \\
         &   2 &  28 & 1.000000e-11 & 5.915080e-10 \\
         &   3 &  32 & 1.000000e-12 & 4.878532e-11 \\
\hline
    1000 &   1 &   0 & 1.000000e-10 & 3.899208e-02 \\
         &   2 &  28 & 1.000000e-11 & 5.915100e-10 \\
         &   3 &  32 & 1.000000e-12 & 4.878548e-11 \\
\hline
   10000 &   1 &   0 & 1.000000e-10 & 3.899208e-02 \\
         &   2 &  28 & 1.000000e-11 & 5.915100e-10 \\
         &   3 &  32 & 1.000000e-12 & 4.878548e-11 \\
\hline
  100000 &   1 &   0 & 1.000000e-10 & 3.899208e-02 \\
         &   2 &  28 & 1.000000e-11 & 5.915099e-10 \\
         &   3 &  32 & 1.000000e-12 & 4.878540e-11 \\
\hline
 1000000 &   1 &   0 & 1.000000e-10 & 3.899208e-02 \\
         &   2 &  28 & 1.000000e-11 & 5.914422e-10 \\
         &   3 &  32 & 1.000000e-12 & 4.871775e-11 \\
\hline
\end{tabular}
\caption{\label{TabM} Convergence of the algorithm for $\tau_{gl}=10^{-10}$
and different values of $M$.}
\end{table}

\begin{comment}
MatLab commands for this table:
lnl(  100,0, 1e-10);
lnl( 1000,0, 1e-10);
lnl(10000,0, 1e-10);
...
\end{comment}

Table~\ref{TabM} shows how the algorithm given above performs.  After 3
iterations, the desired accuracy $\tau_{gl}=10^{-10}$ is achieved.  Moreover, we
can see from the table that the number of iterations are independent of $M$,
that means the algorithm behaves robustly with respect to the problem size $M$.

\begin{figure}
\includegraphics[width=0.5\textwidth]{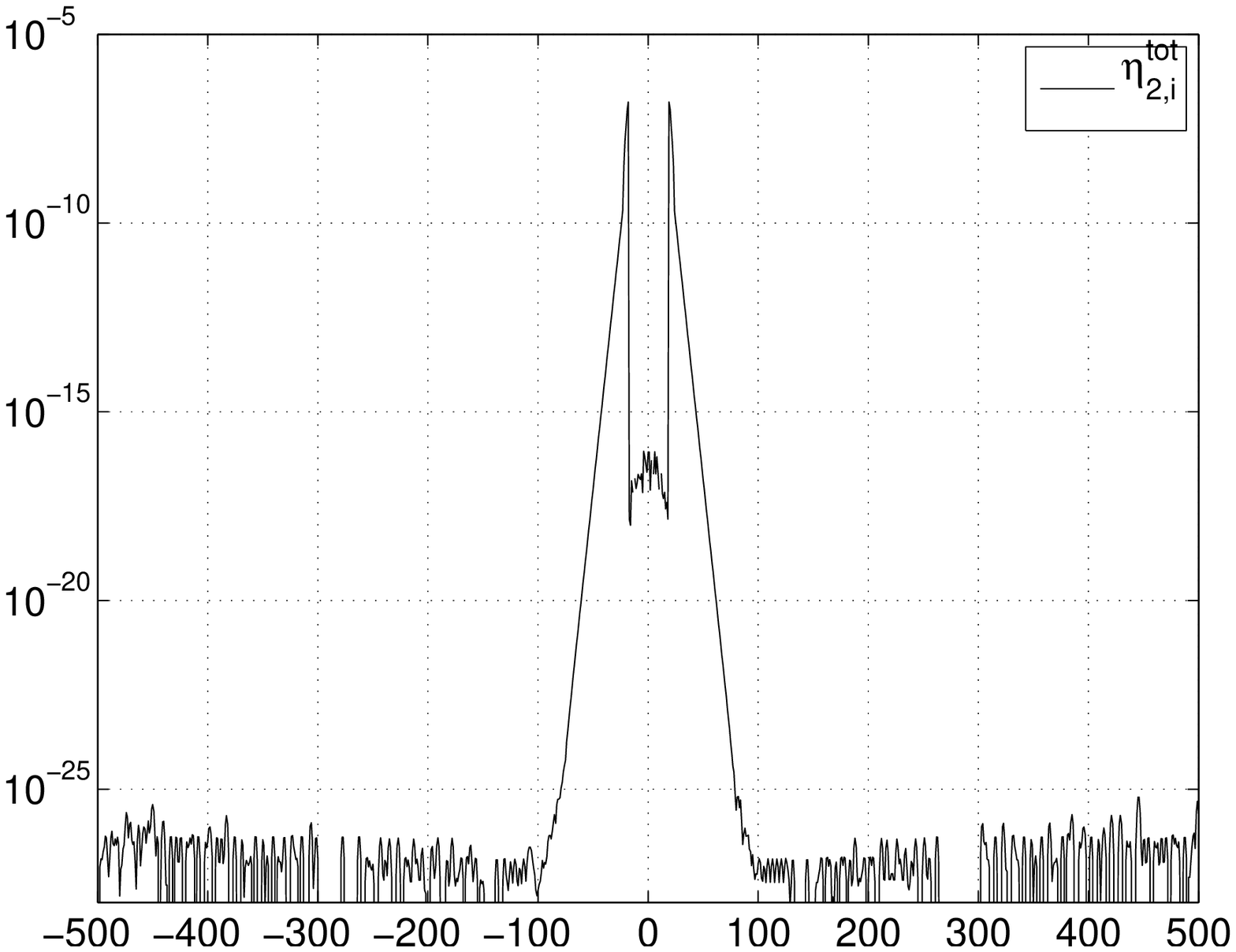}%
\includegraphics[width=0.5\textwidth]{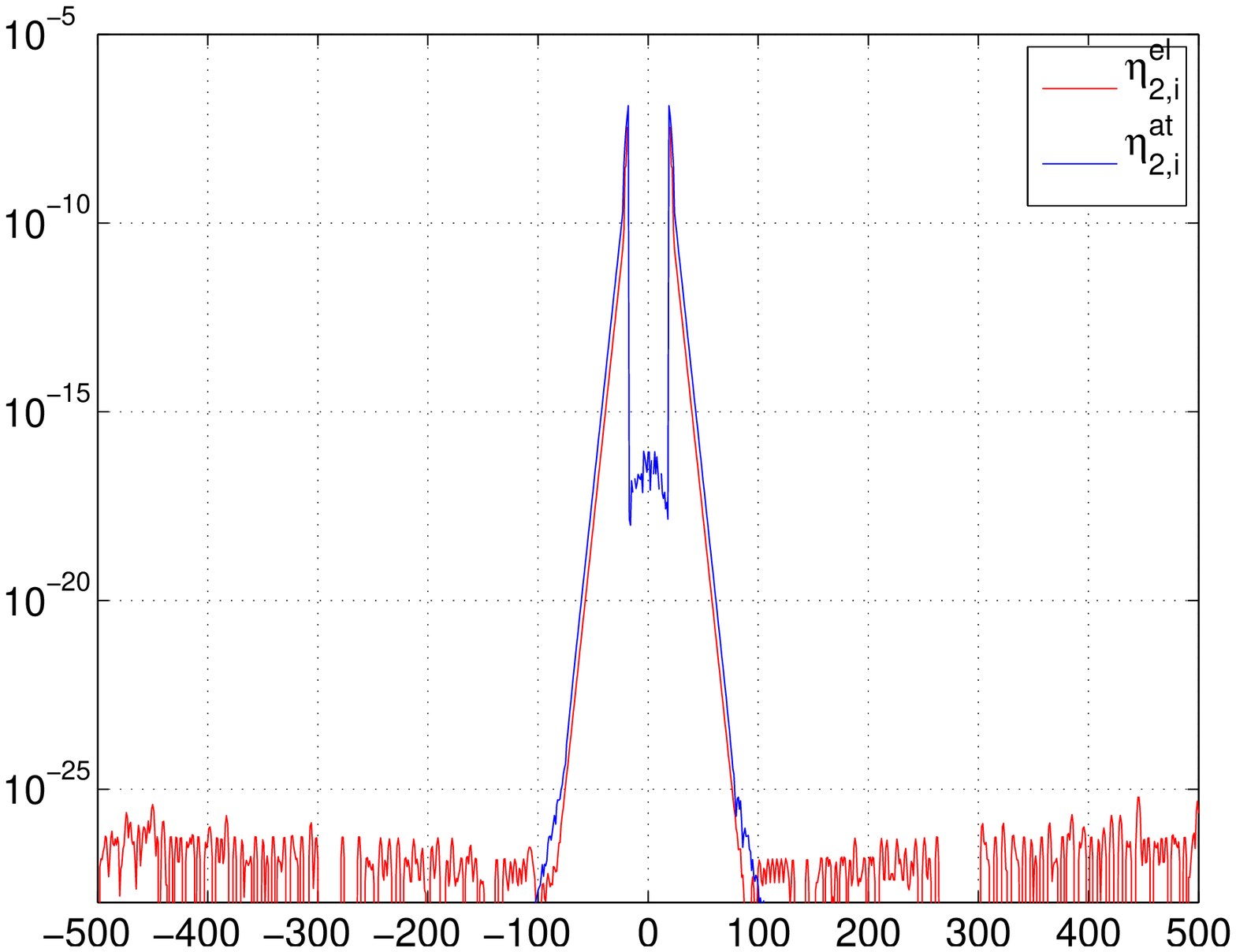}
\caption{The error estimators $\eta_{2,i}^{tot}$ (left) and
$\eta_{2,i}^{el}, \eta_{2,i}^{at}$ (right) for $M=500$, $K=20$.}
\label{FigAtomWise}
\end{figure}

\begin{comment}
MatLab commands:

set(gca, 'fontsize', 16);
[eta_el2, eta_at2, eff] = lnl(500, 20, 1e10);
semilogy(-499:500, eta_el2+eta_at2, 'k');
axis([-500 500 1e-28 1e-5])
set(gca, 'xtick', -500:100:500);
legend('\eta^{tot}_{2,i}');
grid on
print -depsc2 err_div.eps
semilogy(-499:500, eta_el2, 'r', -499:500, eta_at2, 'b');
axis([-500 500 1e-28 1e-5])
set(gca, 'xtick', -500:100:500);
legend('\eta^{el}_{2,i}', '\eta^{at}_{2,i}');
grid on
print -depsc2 err_div_elat.eps
\end{comment}

Figure~\ref{FigAtomWise} (left) shows the decomposition of the error estimator $\eta_2$ for a
typical setting $M=500$, $K=20$. One can clearly see that the error in the
atomistic region is small, whereas the error is large in the
continuum regions that border the atomistic region.
It then decreases exponentially towards the endpoints.

The error in both the atomistic region around the center and the continuum
regions far away from the center are in the range of the (relative) machine
precision $\eps_{mach}$, which accounts for the fluctuations in these
regions. The error can be considered to be numerically zero in these regions. In
the continuum regions, we observe an error of magnitude $\Oh(\eps_{mach}^2)$,
whereas in the continuum region we have $\Oh(\eps_{mach})$, which leads to the
different magnitudes of the fluctuations.

Figure~\ref{FigAtomWise} (right) shows the element-wise contributions
$\eta_{2,i}^{el}$ and the atom-wise contributions $\eta_{2,i}^{at}$ of the
decomposed error estimator $\eta_{2,i}^{tot}=\eta_{2,i}^{at} + \half \big(
\eta_{2,i-1}^{el} + \eta_{2,i}^{el} \big)$. The atomistic part
$\eta_{2,i}^{at}$, which corresponds to the $\vec g^{acT} R^a(\vec y^{ac})$
term, is dominant in the sense that it is about ten times larger than
$\eta_{2,i}^{el}$, which comes from the estimate for the perturbation term
$\hat{\vec e}^T M^a \vec e$.  The fluctuations due to the limited machine
precion in the atomistic region come from $\eta_{2,i}^{at}$, whereas those in
the continuum region away from the defect stem from $\eta_{2,i}^{el}$.  Let us
note that in other applications of duality-based error estimation, the first term might
not always be the dominant term. For example, in mesh refinement for classical
linear finite elements, the first term even vanishes due to Galerkin
orthogonality.

\begin{table}
\begin{tabular}{|c|c|c|c|c|c|}
\hline
 K & $|Q(\vec y^a-\vec y^{ac})|$ & $\eta_1$ & $\eta_1/|Q(\vec y^a-\vec y^{ac})|$
                              & $\eta_2$ & $\eta_2/|Q(\vec y^a-\vec y^{ac})|$ \\
\hline
 0 & 3.627633e-02 & 3.899208e-02 & 1.074863 & 3.999783e-02 & 1.102588 \\
 2 & 3.375762e-02 & 3.872272e-02 & 1.147081 & 5.101700e-02 & 1.511274 \\
 4 & 3.468605e-03 & 4.343595e-03 & 1.252260 & 5.422007e-03 & 1.563166 \\
 6 & 5.418585e-04 & 7.156249e-04 & 1.320686 & 9.187940e-04 & 1.695635 \\
 8 & 1.227067e-04 & 1.675383e-04 & 1.365356 & 2.193196e-04 & 1.787348 \\
10 & 3.287188e-05 & 4.540984e-05 & 1.381419 & 5.984186e-05 & 1.820457 \\
15 & 1.416914e-06 & 1.966114e-06 & 1.387603 & 2.597488e-06 & 1.833201 \\
20 & 6.267636e-08 & 8.695824e-08 & 1.387417 & 1.148736e-07 & 1.832805 \\
25 & 2.770161e-09 & 3.843388e-09 & 1.387424 & 5.077204e-09 & 1.832819 \\
30 & 1.224369e-10 & 1.698739e-10 & 1.387440 & 2.244073e-10 & 1.832840 \\
35 & 5.410783e-12 & 7.508365e-12 & 1.387667 & 9.918687e-12 & 1.833133 \\
40 & 2.379208e-13 & 3.318024e-13 & 1.394592 & 4.383361e-13 & 1.842362 \\
45 & 8.992806e-15 & 1.430733e-14 & 1.590975 & 1.901601e-14 & 2.114580 \\
50 & 7.771561e-16 & 4.120094e-16 & 0.530150 & 6.201285e-16 & 0.797946 \\
\hline
\end{tabular}
\caption{\label{TabEfficiency} Efficiency of the error estimators, $\eta_1/|Q(\vec y^a-\vec y^{ac})|$
and $\eta_2/|Q(\vec y^a-\vec y^{ac})|,$ for $M=1000$. For $K=45$ and $K=50$ the results
become inaccurate due to limited machine precision.}
\end{table}

\begin{figure}
\includegraphics[width=0.8\textwidth]{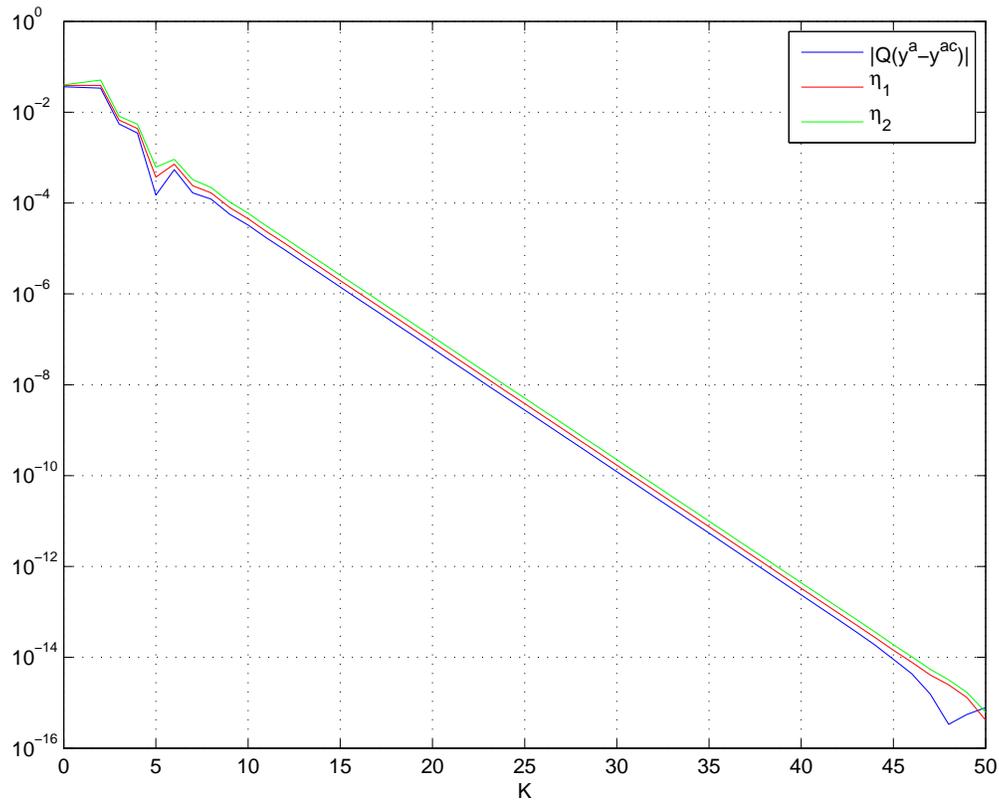}
\caption{Efficiency of the error estimators for $M=1000$.}
\label{FigEfficiency}
\end{figure}

\begin{comment}
MatLab commands:

allK = [0 2:50];            % values for figure
% allK = [0:2:10 15:5:50];  % values for table
efficiency = [];
for K=allK;
    [eta_el2, eta_at2, eff_tmp] = lnl(1000, K, 1e10);
    efficiency = [efficiency; eff_tmp];
end
semilogy(allK, abs(efficiency(:,1)), 'b', ...
         allK, abs(efficiency(:,2)), 'r', ...
         allK, abs(efficiency(:,3)), 'g' );
legend('|Q(y^a-y^{ac})|','\eta_1','\eta_2');
xlabel('K');
grid on;
print -depsc2 efficiency.eps
\end{comment}

Table~\ref{TabEfficiency} and Figure~\ref{FigEfficiency}, which display the same
data, show the efficiency of the error estimators, $\eta_1/|Q(\vec y^a-\vec y^{ac})|$
and $\eta_2/|Q(\vec y^a-\vec y^{ac})|,$
for $M=1000$. For comparison,
the actual error is given as well. For the relatively small 1D problem, the
actual error can be easily computed, whereas in real world applications it is of
course not available.  One can clearly see that $\eta_1$ gives a better estimate
than $\eta_2$, which numerically confirms our conjecture that $\eta_1$ is a
better estimator than $\eta_2$. We see that $\eta_1$ overestimates the actual error by a
factor of 1.4, while $\eta_{2}$ is in a still acceptable range of 2 times the
actual error. Moreover, we can see from Table~\ref{TabEfficiency} and
Figure~\ref{FigEfficiency} that the error decreases exponentially with $K$.

\begin{table}
\begin{tabular}{|c||c|c|c|c|}
\hline
$\tau_{gl}$ & optimal $K$ & $K$ by $\eta_1$ & $K$ by $\eta_2$ \\
\hline
1e-02 &  3 &  3 &  3 \\
1e-03 &  5 &  5 &  5 \\
1e-04 &  9 &  9 & 10 \\
1e-05 & 12 & 13 & 13 \\
1e-06 & 16 & 17 & 17 \\
1e-07 & 20 & 20 & 21 \\
1e-08 & 23 & 24 & 24 \\
1e-09 & 27 & 28 & 28 \\
1e-10 & 31 & 31 & 32 \\
1e-11 & 35 & 35 & 35 \\
1e-12 & 38 & 39 & 39 \\
1e-13 & 42 & 42 & 43 \\
1e-14 & 45 & 46 & 47 \\
\hline
\end{tabular}
\caption{\label{TabOptK} Efficiency of the error estimators for $M=1000$.}
\end{table}

\begin{comment}
MatLab commands:
for K=2:50; lnl(1000, K, 1e10); end
then see when respective error drops below respective tolerance
\end{comment}

Finally, we compare the optimal (smallest) value of $K$ which is needed to
achieve a given accuracy $\tau_{gl}$ with the values for $K$ determined by the
error estimators $\eta_1$ and $\eta_{2}$, again taking into
account the precise error which is available for the model problem.  We see from
Table~\ref{TabOptK} that
even $\eta_2$ only overestimates $K$ by at most 2 atoms. Thus, we get an
efficient estimate of the required atomistic region for both error estimators.

%% file: matrix.tex
% $Id: matrix.tex,v 1.6 2007/03/28 20:16:35 luskin Exp $

%%%%%%%%%%%%%%%%%%%%%%%%%%%%%%%%%%%%%%%%%%%%%%%%%%%%%%%%%%%%%%%%%%%%%%%%%%%%%%%%
%%%%%%%%%%%%%%%%%%%%%%%%%%%%%%%%%%%%%%%%%%%%%%%%%%%%%%%%%%%%%%%%%%%%%%%%%%%%%%%%
%%%%%%%%%%%%%%%%%%%%%%%%%%%%%%%%%%%%%%%%%%%%%%%%%%%%%%%%%%%%%%%%%%%%%%%%%%%%%%%%
\section{Matrix Definitions} \label{AppMatrixDef}

We describe the matrices from Section~\ref{SecMatrixFormulation}.
The matrix
\begin{align}
  D^a = \begin{bmatrix}
    -1 &  1 \\
       & -1 & 1 \\
       &    & \ddots & \ddots \\
       &    &        & -1 & 1
  \end{bmatrix} \in \R^{(2M-1) \times 2M}
\end{align}
transforms atomistic positions to distances between adjacent atoms.  Similarly,
\begin{align}
  D^{qc} = \begin{bmatrix}
     -\nu_{-N+1}^{-1} & \nu_{-N+1}^{-1} \\
    & -\nu_{-N+2}^{-1} & \nu_{-N+2}^{-1} \\
    && \ddots & \ddots \\
    &&& -\nu_{N-1}^{-1} & \nu_{N-1}^{-1}
  \end{bmatrix} \in \R^{(2N-1) \times 2N}
\end{align}
transforms repatom positions from a coarsened chain to normalized distances
between adjacent repatoms.  The matrices
\begin{align}
  (E^a)_{ij} = \begin{cases}
    k_1 +  k_2       & i = j \in \{-M+1,M-1\} \\
    k_1 + 2k_2       & i = j \in \{-M+2,\ldots,M-2\} \\
           k_2       & |j-i|=1 \\
    0 & \text{otherwise},
  \end{cases}
\end{align}
\begin{align}
  (E^{ac})_{ij} = \begin{cases}
    \half k_{12} (\delta_i^c + \delta_{i+1}^c) +
    \half k_1 (\delta_i^a + \delta_{i+1}^a) +
    \half k_2 (\delta_{i-1}^a+\delta_{i}^a+\delta_{i+1}^a+\delta_{i+2}^a)
                     & i = j \\
    \half k_2 (\delta_i^a    + \delta_{i+2}^a)   & j=i+1 \\
    \half k_2 (\delta_{i-1}^a + \delta_{i+1}^a)   & j=i-1 \\
    0 & \text{otherwise},
  \end{cases}
\end{align}
and
\begin{align}
  (E^{qc})_{ij} = \begin{cases}
    \omega_i k_{12} + \half k_1 (\delta_i^a + \delta_{i+1}^a)
       + \half k_2 (\delta_{i-1}^a+\delta_{i}^a+\delta_{i+1}^a+\delta_{i+2}^a)
                     & i = j \\
    \half k_2 (\delta_i^a    + \delta_{i+2}^a)   & j=i+1 \\
    \half k_2 (\delta_{i-1}^a + \delta_{i+1}^a)   & j=i-1 \\
    0 & \text{otherwise},
  \end{cases}
\end{align}
for $i,j=-M+1,\ldots, M-1$ and $i,j = -N+1, \ldots, N-1$, respectively, describe
the spring interactions in terms of the distances between atoms or repatoms.
Accordingly, the matrices
\begin{align}
  K^a = \begin{bmatrix}
    k_0 \\
    & \ddots \\
    && k_0
\end{bmatrix} \in \R^{2M \times 2M}
\end{align}
and
\begin{align}
  (K^{qc})_{ij} = \begin{cases}
    \frac{1}{6} k_0 \left[ (2\nu_{i-1}+\nu_{i-1}^{-1})
                         + (2\nu_{i  }+\nu_{i  }^{-1}) \right] & i = j \in \{ -N+2,\ldots,N-1\}\\
    \frac{1}{6} k_0 (2\nu_{-N+1}+\nu_{-N+1}^{-1}) & i = j = -N+1 \\
    \frac{1}{6} k_0 (2\nu_{N-1}+\nu_{N-1}^{-1})   & i = j = N \\
    \frac{1}{6} k_0 (\nu_i-\nu_i^{-1})           & j = i+1 \\
    \frac{1}{6} k_0 (\nu_j-\nu_j^{-1})           & j = i-1 \\
    0                                           & \text{otherwise},
  \end{cases}
\end{align}
for $i,j = -N+1, \ldots, N$ describe the misfit interactions for the
original atomistic system and the QC approximation.  Finally, the constant vectors
\begin{subequations}
\begin{align}
  \vec a^a = & \begin{bmatrix}
    (-M+1) a_0 & (-M+2)a_0 & \cdots & (M-1) a_0 & M a_0
  \end{bmatrix}^T \in \R^{2M}, \\
  \vec a^{qc} = & \begin{bmatrix}
    \ell_{-N+1} a_0 & \ell_{-N+2} a_0 & \cdots & \ell_{N-1} a_0 & \ell_N a_0
  \end{bmatrix}^T \in \R^{2N}, \\
  \vec b^a = & \begin{bmatrix}
    -M a_0 & (-M+1) a_0 & \cdots & -a_0 & a_0 & \cdots & (M-1) a_0 & M a_0
  \end{bmatrix}^T \in \R^{2M}, \\
  \vec b^{qc} = & \begin{bmatrix}
    (\ell_{-N+1}-1) a_0 & (\ell_{-N+2}-1) a_0 & \cdots & (\ell_0-1) a_0 & \ell_1 a_0 & \cdots & \ell_{N-1} a_0 & \ell_N a_0
  \end{bmatrix}^T \in \R^{2M},
\end{align}
\end{subequations}
fix the equilibrium positions for the spring interactions and the misfit
energy, respectively.